\RequirePackage{fix-cm}
\documentclass[smallextended]{svjour3}       
\smartqed  
\usepackage{dpr,fec}
\usepackage[noend]{algpseudocode}
\usepackage{xcolor}
\hypersetup{colorlinks}
\usepackage[justification=centering]{caption}

\def\abovestrut#1{\rule[0in]{0in}{#1}\ignorespaces}

\def\abovespace{\abovestrut{0.12in}}

\newcommand{\AlgBreak}{\vspace*{-.7\baselineskip}\Statex\hspace*{\dimexpr-\algorithmicindent-2pt\relax}\rule{\textwidth}{0.4pt}}

\algrenewcommand{\algorithmiccomment}[1]{\hfill[#1]}

\usepackage{graphicx} 

\usepackage{graphicx}
\usepackage{subcaption}
\usepackage{float}
\begin{document}

\title{Proximal Quasi-Newton Methods for Regularized Convex Optimization with Linear and Accelerated Sublinear Convergence Rates
}

\titlerunning{Proximal Quasi-Newton Methods for Regularized Convex Optimization}        

\author{Hiva Ghanbari         \and
        Katya Scheinberg 
}


\institute{Katya Scheinberg \at
              Dept. of Industrial and Systems Engineering, Lehigh University, Bethlehem, PA, USA. \\
                           \email{katyascheinberg@gmail.com}         \\  
                           The work of this author is partially supported by NSF Grants DMS 13-19356, CCF-1320137,
 AFOSR Grant FA9550-11-1-0239, and  DARPA grant FA 9550-12-1-0406 negotiated by AFOSR. 
           \and
           Hiva Ghanbari \at
              Dept. of Industrial and Systems Engineering, Lehigh University, Bethlehem, PA, USA. \\
               \email{hiva.ghanbari@gmail.com}
}

\date{October 17, 2017}

\maketitle

\begin{abstract}
A general, inexact, efficient proximal quasi-Newton algorithm for composite optimization problems has been proposed by Scheinberg and Tang [\textit{Math. Program.,} 160~(2016), pp.~495-529] and a sublinear global convergence rate has been established. In this paper, we analyze the global convergence rate of this method, both in the exact and inexact setting, in the case when the objective function is strongly convex. We also investigate a practical variant of this method by establishing a simple stopping criterion for the subproblem optimization. Furthermore, we consider an accelerated variant, based on FISTA of Beck and Teboulle [\textit{SIAM J. Optim.,} 2~(2009), pp.~183-202], to the proximal quasi-Newton algorithm. Jiang, Sun, and Toh [\textit{SIAM J. Optim.,} 22~(2012), pp.~1042-1064] considered a similar accelerated method, where the convergence rate analysis relies on very strong impractical assumptions on Hessian estimates. We present a modified analysis while relaxing these assumptions and perform a numerical comparison of the accelerated proximal quasi-Newton algorithm and the regular one. Our analysis and computational results show that acceleration may not bring any benefit in the quasi-Newton setting.

\end{abstract}



\section{Introduction}
We address the convex optimization problem of the form
\bequation\label{f0}
\min_x \{F(x) := f(x)+g(x),~ x\in \mathbb{R}^n\},
\eequation
where $g:\mathbb{R}^n \rightarrow \mathbb{R}$ is a continuous convex function which is possibly nonsmooth and $f:\mathbb{R}^n \rightarrow \mathbb{R}$ is a convex continuously differentiable function with Lipschitz continuous gradient, i.e.
\bequationn
\|\nabla f(x) -\nabla f(y)\|_2 \leq L \|x-y\|_2~~~\text{for every}~~x,y\in \mathbb{R}^n,
\eequationn
where $L$ is the global Lipschitz constant of the gradient $\nabla f$. This class of problems, when $g(x) = \lambda \|x\|_1$, contains some of the most common machine learning models such as sparse logistic regression \cite{lin,shalev}, sparse inverse covariance selection \cite{hsieh,olsen,rish}, and unconstrained Lasso \cite{tibshirani}. 

The {Proximal Gradient Algorithm (PGA)} is a variant of the \textit{proximal methods} and is a well-known first-order method for solving optimization problem \eqref{f0}. Although  classical \textit{subgradient methods} can be applied to problem \eqref{f0} when $g$ is nonsmooth, they can only achieve the rate of convergence of $\mathcal{O}(1/\sqrt{k})$ \cite{nesterov2}, while {PGA} converges at a rate of $\mathcal{O}(1/k)$ in both smooth and nonsmooth cases \cite{nesterov4,beck}. In order to improve the global sublinear rate of convergence of {PGA} further, the {Accelerated Proximal Gradient Algorithm (APGA)} has been originally proposed by Nesterov in \cite{nesterov3}, and refined by Beck and Teboulle in \cite{beck}. It has been shown that the {APGA}  provides a significant improvement compared to {PGA}, both theoretically, with a rate of convergence of $\mathcal{O}({1}/{k^2})$,  and practically \cite{beck}. This rate of convergence is known to be the best that can be obtained using only first-order information \cite{yudin,nesterov2,goldfarb}, causing {APGA} to be known as an \textit{optimal first-order method}. The  class of accelerated methods contains many variants that share the same convergence rates and use only first-order information \cite{nesterov2,nesterov5,tseng}. The main known drawback of most of the variants of {APGA} is that the sequence of the step-size $\{\mu_k\}$ has to be nonincreasing. This theoretical restriction sometimes has a big impact on the performance of this algorithm in practice. In \cite{goldfarb}, in order to overcome this difficulty, a new version of {APGA} has been proposed. This variant of  {APGA} allows to increase step-sizes from one iteration to the next, but maintain  the same rate of convergence of $\mathcal{O}({1}/{k^2})$. In particular, the authors have shown that a full backtracking strategy can be applied in {APGA} and that the resulting complexity of the algorithm depends on the average value of step-size parameters, which is closely related to local Lipschitz constants, rather than the global one.
  
To make {PGA} and {APGA} practical, for some complicated instances of \eqref{f0}, one needs to allow for inexact computations in  the algorithmic steps.
In \cite{mark}, inexact variants of {PGA} and {APGA} have been analyzed with two possible sources of error: one in the calculation of the gradient of the smooth term and the other in the proximal operator with respect to the nonsmooth part. The convergence rates are preserved if
the sequence of errors converges to zero sufficiently fast.  Moreover, it has been shown that both of these algorithms obtain a linear rate of convergence, 
when the smooth term $f$ is strongly convex\footnote{For APGA a different variant is analyzed in the case of strong convexity.}. Recently, in \cite{dmitri}, the linear convergence of {PGA} has been shown under the \textit{quadratic growth condition}, which is weaker than a strong convexity assumption. In particular, their analysis relies on the fact that {PGA} linearly bounds the distance to the solution set by the step lengths. This property, called an \textit{error bound condition}, has been proved to be equivalent to the {standard quadratic growth condition}. More precisely, a strong convexity assumption is a sufficient, but not a necessary condition for this error bound property.

 While {PGA} and {APGA} can be efficient in solving \eqref{f0}, it has been observed that using (partial) second-order information often significantly improves the performance of the algorithms. Hence,  Newton type  proximal algorithms, also known as the \textit{proximal Newton methods}, have become popular  \cite{sra,olsen,lee,nocedal3,tang} and are often the methods of choice. When accurate (or nearly accurate) second-order information is used, the method no longer falls in the first-order category and faster convergence rates are expected, at least locally. Indeed, the global convergence and the local superlinear rate of convergence of the {Proximal Quasi-Newton Algorithm (PQNA)} are presented in  \cite{lee} and \cite{nocedal3}, respectively for both the exact and inexact settings. 
 However, in the case of limited memory BFGS method \cite{nocedal1,nocedal2}, the method is still essentially a first-order method. While practical performance may be by far superior,
 the rates of convergence are at best the same as those of the pure first-order counterparts. In \cite{tang}, an inexact {PQNA} with global sublinear rate 
  of $\mathcal{O}({1}/{k})$ is proposed. While the algorithm can use any positive definite Hessian estimates, as long as their eigenvalues are uniformly bounded above and away from zero, the practical implementation proposed in \cite{tang} used a limited memory BFGS Hessian approximation. The inexact setting of the algorithm
  allows for a relaxed sufficient decrease condition as well as an inexact  subproblem optimization, for example via  coordinate descent. 
 
In this work, we show that {PQNA}, as proposed in \cite{tang}, using general Hessian estimates $H_k$, has the linear convergence rate in the case of strongly convex smooth term $f$.  Moreover,  we consider an inexact variant, similar to the ones in \cite{nocedal3,tang}, allowing inexact subproblem solutions as well as a relaxed  sufficient decrease condition. In order to control the errors in the inexact subproblem optimization, we establish a simple stopping criterion for the subproblem solver, based on the iteration count, which guarantees that the inner subproblem is solved to the required  accuracy. In contrast, in related works \cite{toh,villa}, it is assumed that an approximate subproblem solver yields an approximate subdifferential, which is a strong assumption on the subproblem solver which also does not clearly result in a simple  stopping criterion. 
  
Next, we apply  Nesterov's acceleration scheme to {PQNA} as proposed in \cite{tang}, with a view of developing a version of this algorithm with faster convergence rates in the general convex case. In \cite{toh}, the authors have introduced the {Accelerated Proximal Quasi-Newton Algorithm (APQNA)} with rate of convergence of $\mathcal{O}({1}/{k^2})$. However, this rate of convergence is achieved under condition $0 \prec H_k \preceq H_{k-1}$, on the Hessian estimate $H_k$, at each iteration $k$. At the same time, this sequence of matrices has to be chosen so that $H_k$ is sufficiently positive definite to provide an
overapproximation of $f$.  Hence, these two conditions may contradict with each other unless the  sequence of $\{H_k\}$ consists of unnecessarily large matrices.  Moreover, in a particular case, when $H_k$ is set to be a scalar multiple of the identity, i.e., $H_k = \frac{1}{\mu_k} I$, then assumption $0 \prec H_k \preceq H_{k-1}$ enforces $\mu_k\geq \mu_{k-1}$, implying nondecreasing step-size parameters,  which contradicts the standard condition of {APGA}, which is $\mu_k\leq \mu_{k-1}$.
 
In this work, we introduce a new variant of {APQNA}, where we relax the restrictive assumptions imposed in  \cite{toh}.
We use the scheme, originally introduced in \cite{goldfarb}, which allows for the increasing and decreasing step-size parameters. 
We show that  our version of {APQNA} achieves the convergence rate of $\mathcal{O}(1/k^2)$ under some assumptions on the Hessian estimates. While we show that this assumption is rather strong and may not be satisfied by general matrices, it is not contradictory. 
Firstly, our result  applies under the same restrictive condition from \cite{toh}.
We also show that our condition on the matrices holds  in the case when the approximate Hessian at each iteration is a scaled version of the same ``fixed" matrix $H$, which is a generalization of {APGA}.  We  investigate the performance of this algorithm in practice, and discover that this restricted version
of a proximal quasi-Newton method is quite effective in practice. We also demonstrate that the general L-BFGS based {PQNA} does not benefit from the acceleration, which supports our analysis of the theoretical limitations. 

This paper is organized as follows. In the next section, we describe the basic definitions  and existing algorithms, {PGA}, {APGA} and {PQNA},  that we refer to later in the paper. In Section \ref{ch_strong}, we analyze {PQNA} in the strongly convex case. In Section \ref{ch_accelerated}, we propose, state and analyze a general {APQNA} and its simplified version. We present computational results in Section \ref{sec.result}.  Finally, we state our conclusions in Section \ref{sec.conclusion}.

\section{Notation and Preliminaries} \label{1.preliminaries}

In this work, the Euclidean norm $\|a\|^2 := a^T a$, and the inner product $\langle a, b \rangle := a^T b$, are also defined in the scaled setting such that, $\|a\|^2_H := a^T  H a$, and $\langle a, b \rangle_H := a^T H b$. We denote the identity matrix by $I \in \mathbb{R}^{n \times n}$. The vector $e_j$ stands for a unit vector along the $j$-th coordinate. We use $x_k$ to denote the approximate minimizer (the iterate), computed  at iteration $k$ of an appropriate algorithm, and $x_*$ to denote an exact minimizer of $F$. Finally, $(\partial F(x))_{\text{min}}$ denotes the minimum norm subgradient of function $F$ at point $x$.

The \textit{proximal mapping} of a convex function $g$ at a given point $v$, with parameter $\mu$ is defined as
\bequation \label{defprox}
prox_g^{\mu}(v) := \arg \min_{u \in \mathbb{R}^n} \{g(u) + \frac{1}{2 \mu} \|u-v\|^2\},~~~\text{where}~~ \mu >0.
\eequation
The proximal mapping is the base operation of the \textit{proximal methods}. In order to solve the composite problem \eqref{f0}, each iteration of {PGA} computes the proximal mapping of the function $g$ at point  $x_k- \mu_k \nabla f(x_k)$ as  follows:
\bequation \label{Prox_0}
\begin{aligned}
p_{\mu_k}(x_{k}) &:= prox_g^{\mu_k}(x_k - \mu_k \nabla f(x_k))\\
&:= \arg \min_{u \in \mathbb{R}^n} \{ g(u) + f(x_k) +\langle \nabla f(x_k) ,u-x_k \rangle + \frac{1}{2 \mu_k} \|u - x_k\|^2 \}.
\end{aligned}
\eequation
 We will call the objective function, that is minimized in  \eqref{Prox_0}, a \textit{composite quadratic approximation} of the convex function $F(x) :=f(x)+g(x)$. This approximation at a given point $v$, for a given $\mu$  is defined as
\bequation \label{Q.mu}
Q_{\mu}(u,v) := f(v) + \nabla f(v)^T (u-v)+\frac{1}{2\mu} \|u-v\|^2+g(u).
\eequation
 Then, the proximal operator can be written as 
\bequationn
p_\mu(v) := \arg \min_{u \in \mathbb{R}^n} Q_\mu(u,v).
\eequationn

Using this notation we first present the basic {PGA} framework  with backtracking over $\mu$ in Algorithm \ref{p.al0}. The simple backtracking scheme enforces that the  sufficient decrease condition 
\bequation \label{suff}
F(x_{k+1}) \leq Q_{\mu_k}(x_{k+1},x_k) \leq Q_{\mu_k}(x_{k},x_k) = F(x_k),
\eequation
holds. This condition is essential in the convergence rate analysis of {PGA} and is easily satisfied when $\mu\leq {1}/{L}$. The backtracking is used for two reasons--because the constant $L$ may not be known and because $\mu\leq {1}/{L}$ may be too pessimistic, i.e., condition \eqref{suff} may be satisfied for much larger values of $\mu$ allowing for larger steps. 


\balgorithm[ht]
  \small
  \caption{Proximal Gradient Algorithm}
  \label{p.al0}
  \balgorithmic[1]
    \smallskip
    \AlgBreak
      \State Initialize $x_0 \in \mathbb{R}^n$, and choose $\beta \in (0,1)$.
      \For{$k=1,2,\cdots$}
      \State Choose $\mu_k^0$ and define $\mu_k := \mu_k^0$.
      \State Compute $p_{\mu_k}(x_k) :=\arg \min_{u\in \mathbb{R}^n} Q_{\mu_k}(u,x_k)$.
      \While{$F\(p_{\mu_k}(x_k)\) > Q_{\mu_k}\(p_{\mu_k}(x_k),x_k\)$}
      \State Set $\mu_k \gets \beta \mu_k$.
      \State Compute $p_{\mu_k}(x_k) :=\arg \min_{u\in \mathbb{R}^n} Q_{\mu_k}(u,x_k)$.
      \EndWhile
      \State Set $x_{k+1} \gets p_{\mu_k}(x_k)$.
      \EndFor
  \ealgorithmic
\ealgorithm


We now present  the accelerated variant of {PGA} stated as {APGA}, where at each iteration $k$, instead of constructing $Q_{\mu_k}$ at the current minimizer $x_k$, it is constructed at a different point $y_k$, which is chosen as a specific linear combination of the latest two or more minimizers, 
e.g. 
\bequationn \label{acc.gen}
y_{k+1} = x_k + \alpha_k (x_k - x_{k-1}),
\eequationn
where the sequence $\{ \alpha_k \}$ is chosen in such a way to guarantee an accelerated convergence rate as compared to the original {PGA}. 
Algorithm \ref{p.al1} is a variant of {APGA} framework, often referred to as FISTA, presented in \cite{beck}, where $\alpha_k = ({t_{k}-1}) /  ({t_{k+1}})$. In this work, we choose to focus on FISTA algorithm specifically.  The choice of the accelerated parameter $t_{k+1}$ in \eqref{tal2} is dictated by the analysis of the complexity of FISTA \cite{beck} and the condition  $\mu_{k+1}\leq \mu_{k}$ that is imposed by the initialization of the backtracking procedure with $\mu_{k+1}^0:=\mu_k$. In \cite{goldfarb}, the definition of $t_{k+1}$ was generalized to allow $\mu_{k+1}^0>\mu_k$, while retaining the convergence rate. We will use a similar technique in our proposed {APQNA}. 
  

\balgorithm[ht]
  \small
  \caption{Accelerated Proximal Gradient Algorithm}
  \label{p.al1}
  \balgorithmic[1]
    \smallskip
    \AlgBreak
      \State  Initialize $t_1=1$, $\mu_1^0 > 0$, and $y_1=x_0 \in \mathbb{R}^n$, and choose $\beta \in (0,1)$.
      \For{$k=1,2,\cdots$}
      \State Define $\mu_k := \mu_k^0$.
      \State Compute $p_{\mu_k}(y_k) :=\arg \min_{u\in \mathbb{R}^n} Q_{\mu_k}(u,y_k)$.
      \While{$F\(p_{\mu_k}(y_k)\) > Q_{\mu_k}\(p_{\mu_k}(y_k),y_k\)$}
      \State Set $\mu_k \gets \beta \mu_k$.
      \State Compute $p_{\mu_k}(y_k) :=\arg \min_{u\in \mathbb{R}^n} Q_{\mu_k}(u,y_k)$.
      \EndWhile
      \State Set $x_k \gets p_{\mu_k}(y_k)$.
      \State Define $\mu_{k+1}^0:=\mu_k$ and compute $t_{k+1}$ and $y_{k+1}$, so that
\bsubequations 
  \begin{align}
&t_{k+1}= \frac{1}{2} \(1+\sqrt{1+4t_{k}^2}\) \label{tal2}\\
\text{and}~~&y_{k+1}=x_{k}+\frac{t_{k}-1}{t_{k+1}}\(x_{k}-x_{k-1}\). \label{tall2}
\end{align}
\esubequations

       \EndFor
  \ealgorithmic
\ealgorithm


In this work, we are interested in the extensions of the above proximal methods,  which utilize an approximation function $Q_{\mu}$, using partial second-order information about $f$. These quasi-Newton type proximal algorithms use a generalized form of the proximal operator \eqref{defprox}, known as the \textit{scaled proximal mapping} of $g$, which are defined for a given point $v$  as
\bequationn \label{scalprox}
prox_g^H(v) := \arg \min_{u \in \mathbb{R}^n} \{g(u) + \frac{1}{2} \|u-v\|_H^2\},
\eequationn
where matrix $H$ is a positive definite matrix. In particular,  the following operator
\bequation \label{Prox_00}
\begin{aligned}
p_{H_k}(x_k) &:= prox_g^{H_k}\(x_k - H_k^{-1} \nabla f(x_k)\),
\end{aligned}
\eequation
computes the minimizer, over $u$, of the following composite quadratic approximation of function $F$ 
\bequation \label{e.b1}
Q_H(u,v) := f(v) +  \langle \nabla f(v), u-v \rangle + \frac{1}{2} \|u-v\|_H^2 + g(u),
\eequation
when $v=x_k$.  Matrix $H$ is the approximate Hessian of $f$ and its choice plays the key role in the quality of this approximation. Clearly, when $H = \frac{1}{\mu}I$, approximation \eqref{e.b1} converts to \eqref{Q.mu}, which is used throughout {PGA}. If we set $H = \nabla^2 f(v)$, then \eqref{e.b1} is the second-order Taylor approximation of $F$. At each iteration of {PQNA}  the following optimization problem needs to be solved 
\bequation \label{e.b2}
p_H(v) := \arg \min_{u \in \mathbb{R}^n} Q_H(u,v),
\eequation
which we assume to be computationally inexpensive relative to solving \eqref{f0}  for any $v \in \mathbb{R}^n$ and some chosen class of positive definite approximate Hessian $H$. Our assumption is motivated by \cite{tang}, where it is shown that for L-BFGS Hessian approximation, problem \eqref{e.b2} can be solved efficiently  and inexactly via coordinate descent method. Specifically, the proximal operator \eqref{Prox_00} does not have closed form solution for most types of Hessian estimates $H_k$ and most nonsmooth terms $g$, such as $g=\lambda \|x\|_1$, when \eqref{Prox_00} is a convex quadratic optimization problem. Hence, it may be too expensive to seek the exact solution of subproblem \eqref{Prox_00} on every iteration. In \cite{tang}, an efficient version of {PQNA} is proposed that constructs Hessian estimates based on the L-BFGS updates, resulting in $H_k$ matrices that are sum of a diagonal and a low rank matrix. The resulting  subproblem, structured as \eqref{Prox_00}, is then solved up to some expected accuracy via randomized coordinate descent, which effectively exploits the special structure of $H_k$.
 The analysis in \cite{tang} shows that the resulting inexact {PQNA} converges sublinearly if the Hessian estimates remain positive definite and bounded, without assuming any other structure. In this paper, all the theory is derived for arbitrary positive definite Hessian estimates, without any assumption on their structure, or how closely they are representing the true Hessian. In our implementation, however, we will construct the Hessian estimates via L-BFGS as it is done in \cite{tang}. The framework of the inexact variant of {PQNA} for general approximate Hessian is stated in Algorithm \ref{al.pqna}. In this algorithm, 
the inexact solution of  \eqref{e.b2} is denoted by $p_{H,\epsilon}(v)$, as an $\epsilon-$minimizer of the subproblem that satisfies
\bequation \label{i.min}
g(p_{H,\epsilon}(v)) + \frac{1}{2} \| p_{H,\epsilon}(v) - z \|_H^2 \leq \min_{u \in \mathbb{R}^n} \{g(u) + \frac{1}{2} \|u - z \|_H^2\}+ \epsilon,
\eequation
where $z := v - H^{-1} \nabla f(v)$. Obtaining such an inexact solution can be achieved by applying any linearly convergent algorithm, as will be discussing in detail at the end of this section. 

In addition, for a given $\eta \in (0,1]$, the typical condition \eqref{suff}, used in \cite{beck} and \cite{bach}, is relaxed by using a trust-region like sufficient decrease condition
\bequation \label{step.cond}
\( F\(p_{H,\epsilon}(v)\) - F(v) \) \leq \eta \( Q_H\(p_{H,\epsilon}(v),v\) - F(v)\).
\eequation
 This relaxed condition was proposed and tested in \cite{tang} for {PQNA} and was shown to lead to superior numerical performance, saving multiple backtracking steps during the earlier iterations of the algorithm. Note that, one can obtain the exact version of Algorithm \ref{al.pqna} by replacing $p_{H_k,\epsilon_k}$ with $p_{H_k}$, and setting $\eta = 1$.


\balgorithm[ht]
  \small
  \caption{Inexact Proximal Quasi-Newton Algorithm}
 \label{al.pqna}
  \balgorithmic[1]
    \smallskip
    \AlgBreak
      \State  Initialize $x_0 \in \mathbb{R}^n$, and choose $ \beta \in (0,1)$ and $ \eta \in (0,1]$.
      \For{$k=1,2,\cdots$}
      \State Choose $\mu_k > 0$ and bounded $G_k \succeq 0$.
      \State Define $H_k := G_k + \frac{1}{2 \mu_k} I$.
      \State Compute $p_{H_k,\epsilon_k}(x_k)$ such that \eqref{i.min} is satisfied.
      \While{$\(F\(p_{H_k,\epsilon_k}(x_k)\) - F(x_k)\) > \eta \( Q_{H_k}\(p_{H_k,\epsilon_k}(x_k),x_k\) - F(x_k)\)$}
      \State Set $\mu_k \gets {\beta} \mu_k$.
      \State Update $H_k$ via $H_k := G_k + \frac{1}{2 \mu_k} I$.
      \State Compute $p_{H_k,\epsilon_k}(x_k)$ such that \eqref{i.min} is satisfied.
      \EndWhile
      \State Set $x_{k+1} \gets p_{H_k,\epsilon_k}(x_k)$.
       \EndFor
  \ealgorithmic
\ealgorithm


\newpage
Throughout our analysis, we  make the following assumptions.
\bassumption \label{assum.main} 
\ \\
\bitemize
\item  $f$ is convex with Lipschitz continuous gradient with constant $L$.
\item $g$ is a lower semi-continuous proper convex function.
\item There exists an $x_*\in \mathbb{R}^n$, which is a minimizer of $F$.
\item There exists positive constants $m$ and $M$ such that, for all $k >0$,
\bequation \label{bound.strong}
m I \preceq H_k \preceq M I.
\eequation
\eitemize
\eassumption

\bremark  \label{G.bound}
In Algorithm \ref{al.pqna}, as long as the sequence of positive definite matrices $G_k$ has uniformly  bounded eigenvalues, condition \eqref{bound.strong} is satisfied. In fact, since the sufficient decrease condition in Step 3 is satisfied for $H_k \succeq L I$, then it is satisfied when $\mu_k\leq {1}/{L}$. Hence, at each iteration we have a finite and bounded number of backtracking steps and the resulting $H_k$ has bounded eigenvalues. The lower bound on the eigenvalues of $H_k$ is simply imposed either by choosing a positive definite $G_k$ or bounding $\mu_k^0$ from above.
\eremark

In the next section, we analyze the convergence properties of {PQNA} when $f$ in  \eqref{f0} is strongly convex.

\section{Analysis of the Proximal Quasi-Newton Algorithm under Strong Convexity}\label{ch_strong}
In this section, we analyze the convergence properties of {PQNA} to solve problem \eqref{f0}, in the case when the smooth function $f$ is $\gamma-$strongly convex. In particular, the following assumption is made throughout this section.
\bassumption  \label{assum.strong}
For all $x$ and $y$ in $\mathbb{R}^n$, and any $t \in [0,1]$, the following two equivalent conditions hold.
\bsubequations
\begin{align}
& \gamma \|x-y\|^2 \leq \langle  \nabla f(x) - \nabla f(y), x-y  \rangle \label{def.strong1}\\
\text{and}~~&f\(tx+(1-t)y\) \leq t f(x) + (1-t) f(y) - \frac{1}{2} \gamma t(1-t) \|x-y\|^2. \label{def.strong}
\end{align}
\esubequations
\eassumption

To establish a linear convergence rate of {PQNA}, we consider extending two different approaches used to show a similar result for {PGA}. The first approach we consider can be found in \cite{phd}, and is based on the proof techniques used in \cite{dmitri} for {PGA}.  The reason we chose the approach  in \cite{dmitri}  is due to the fact that the linear rate of convergence is shown under the quadratic growth condition, which is a relaxation of the strong convexity. Hence, extending this analysis to {PQNA}, as a subject of a future work, may allow us to  relax the strong convexity assumption for this algorithm as well. However, there appears to be some limitations in the extension of this analysis \cite{phd}, in particular in the inexact case. This observation motivates us to present the  approach used in \cite{nesterov4} to analyze convergence properties of inexact {PQNA}. As we see below, this analysis readily extends to our case and allows us to establish simple rules for subproblem solver termination to achieve the desired subproblem accuracy.

\subsection{Convergence Analysis}
Let us consider Algorithm \ref{al.pqna} for which  \eqref{i.min} holds for some sequence of errors $\epsilon_k \geq 0$. The relaxed sufficient decrease condition 
\bequationn \label{step.cond2}
F(x_{k+1}) - F(x_k) \leq \eta ( Q_{H_k}(x_{k+1},x_k) - F(x_k)),
\eequationn
for a given $\eta \in (0,1]$, can be written as
\bequationn
\baligned
F(x_{k+1}) &\leq Q_{H_k}(x_{k+1},x_k) - (1 - \eta) \( Q_{H_k}(x_{k+1},x_k) - F(x_k) \) \\
& \leq Q_{H_k}(x_{k+1},x_k) -\frac{1-\eta}{\eta} \( F(x_{k+1}) - F(x_k) \).
\ealigned
\eequationn
Thus, at each  iteration we have 
\bequation \label{i.dec}
F(x_{k+1})  \leq Q_{H_k}(x_{k+1},x_k) + \xi_k,
\eequation
where the sequence of the errors $\xi_k $ is defined as
\bequation \label{xi.3}
\xi_k \leq (1 - \frac{1}{\eta})  \( F(x_{k+1}) - F(x_k) \).
\eequation

In particular, setting $\eta=1$ results in $\xi_k=0$, for all $k$ and enforces the algorithm to  accept only those steps that achieve full (predicted) reduction. 
However, using $\eta <1$ allows the algorithm to take steps satisfying only a fraction of the predicted reduction, which may lead to larger steps and faster progress.

Under the above inexact condition, we can show the following result.
\btheorem \label{sec.linear}
Suppose that Assumptions \ref{assum.main} and \ref{assum.strong} hold. At each iteration of the inexact {PQNA}, stated in Algorithm \ref{al.pqna}, we have
\bequation \label{i.rate.strong}
F(x_k) - F(x_*) \leq \rho^k \( F(x_0) - F(x_*) +A_k \),
\eequation
when $\rho = 1 - {(\eta \gamma)}/{(\gamma + M)}$, and 
\bequationn
A_k := \eta \sum_{i=1}^{k} \({\epsilon_i} / {\rho^i} \).
\eequationn

\etheorem
\bproof
 Applying \eqref{i.dec}, with  $v = x_k$ and consequently $p_{H_k,\epsilon_k}(x_k) = x_{k+1}$, we  have
\bequationn
\baligned
F(x_{k+1}) &\leq Q_{H_k}(x_{k+1},x_k) + \xi_k\\
& = f(x_k) + \langle \nabla f(x_k), x_{k+1}-x_k \rangle + \frac{1}{2} \|x_{k+1}-x_k\|_{H_k}^2 + g(x_{k+1})+\xi_k \\
& = \min_{u \in \mathbb{R}^n} ~  f(x_k) +  \langle \nabla f(x_k), u-x_k \rangle + \frac{1}{2} \|u-x_k\|_{H_k}^2 + g(u)+\epsilon_k + \xi_k \\
& \leq  \min_{u \in \mathbb{R}^n} ~ f(u) + \frac{1}{2} \|u-x_k\|_{H_k}^2 + g(u)+(\epsilon_k + \xi_k) ~~~~~\(\text{convexity of $f$}\) \\
& =  \min_{u \in \mathbb{R}^n} ~ F(u) + \frac{1}{2} \|x_k-u\|_{H_k}^2+(\epsilon_k + \xi_k)\\
& \leq \min_{t \in [0,1]} ~ F( t x_* + (1-t) x_k) + \frac{1}{2} \|x_k- t x_* - (1-t) x_k\|_{H_k}^2\\
&~+(\epsilon_k + \xi_k) \\
& \leq \min_{t \in [0,1]} ~ t F(x_*) + (1-t) F(x_k) - \frac{1}{2} \gamma t(1-t) \|x_*-x_k\|^2 \\ &~+ \frac{1}{2} t^2 \|x_*-x_k\|_{H_k}^2+(\epsilon_k + \xi_k) ~~~~(\text{using \eqref{def.strong}}) \\ 
& \leq \min_{t \in [0,1]} ~ t F(x_*) + (1-t) F(x_k) - \frac{1}{2} \gamma t(1-t) \|x_*-x_k\|^2 \\&~ + \frac{1}{2} M t^2\|x_*-x_k\|^2+(\epsilon_k + \xi_k)\\
& \leq  t' F(x_*) + (1-t') F(x_k)+(\epsilon_k + \xi_k). ~~~~~~~~\(\text{where $t' = \frac{\gamma}{\gamma + M}$}\) \\
\ealigned
\eequationn
Therefore, we have
\bequationn
F(x_{k+1}) \leq  t' F(x_*) + (1-t') F(x_k)+(\epsilon_k + \xi_k),
\eequationn
which implies 
\bequationn
\baligned
F(x_{k+1}) - F(x_*) &\leq (1-t') (F(x_k) - F(x_*))+(\epsilon_k + \xi_k).
\ealigned
\eequationn
Now, by substituting the expression for $\xi_k$, as stated in \eqref{xi.3}, we will have
\bequationn \label{ett}
F(x_{k+1}) - F(x_*) \leq \rho (F(x_k) - F(x_*))+ \eta \epsilon_k,
\eequationn
where $\rho = 1-\eta t'$. Now, we can conclude the final result as
\bequationn
\baligned
F(x_k) - F(x_*) &\leq \rho^k(F(x_0) - F(x_*)) + \sum_{i=1}^{k} \eta \rho^{k-i} \epsilon_i  \\
& = \rho^k \(F(x_0) - F(x_*) + \eta \sum_{i=1}^{k} \({\epsilon_i} / {\rho^i} \) \),
\ealigned
\eequationn
where $\rho = 1 - {(\eta \gamma)}/{(\gamma + M)}$.
\eproof

\bremark
In Theorem \ref{sec.linear}, by setting $\epsilon_k = 0$ and $\eta = 1$, which implies $\xi_k = 0$, we  achieve the linear convergence rate of the exact variant of {PQNA}.
\eremark

\bremark
We have shown that the linear rate of {PQNA} is $\rho= 1 - {(\eta \gamma)}/{(\gamma + M)} $. 
As argued in Remark \ref{G.bound}, $M$ is of the same order as $L$ in the worst case, hence in that case the linear rate of {PQNA} is the same as that of the simple {PGA}. However, it is easy to see that in the proof of Theorem \ref{sec.linear},
the linear rate is derived using the upper bound on $\|x_*-x_k\|_{H_k}^2$, where $H_k$ is the approximate Hessian on step $k$. Clearly, the idea of using the partial second-order information is to reduce the worst case bound of $H_k$ in general and consequently on $\|x_*-x_k\|_{H_k}^2$. In particular, obtaining a smaller bound $M_k$ on each iteration yields a larger convergence coefficient $\rho_k=1 - {(\eta \gamma)}/{(\gamma + M_k)}$. While for general $H_k$, we do not expect to improve upon the regular {PGA} in theory, this remark serves to explain the better performance of {PQNA} in practice. 
 \eremark

Based on the result of Theorem \ref{sec.linear}, it follows that the boundedness of the sequence $\{A_k\}$ is a sufficient condition to achieve the linear convergence rate. 
Hence, the required condition on the sequence of errors is  $\sum_{i=1}^{k} \({\epsilon_i} / {\rho^i} \) < \infty$.
For all $i \leq k$, suppose that $ \epsilon_i\leq C \rho^{i \cdot \delta} $, for some $\delta > 1$ and some $C>0$. Then, we have $\sum_{i=1}^{k} \({\epsilon_i} / {\rho^i} \) \leq C \sum_{i=1}^{k} \rho^{i(\delta-1)}$, which is uniformly bounded for all $k$. 
Recall that the $k$-th subproblem $Q_k^* := \min_{u \in \mathbb{R}^n} Q_{H_k}(u,x_k)$ is a strongly convex function with strong convexity parameter at least $m$--the lower bound on the eigenvalues of $H_k$. Let us assume now that each subproblem is solved via an algorithm with a linear convergence rate for strongly convex problems. In particular, if the subproblem solver is applied for $r(k)$ iterations to the $k$-th subproblem, we have
\bequation \label{Q.linear}
\(Q_{H_k}( u_{r(k)},x_k) -Q_k^* \) \leq \alpha \(Q_{H_k}(u_{{r(k)}-1},x_k)-Q_k^*\),
\eequation

where $\alpha \in (0,1)$. Our goal is to ensure that $ \epsilon_k\leq C \rho^{k \cdot \delta} $, which can be achieved by applying sufficient number of iterations of the subproblem algorithm. To be specific, the following theorem characterizes this required bound on the number of inner iterations.
 
\btheorem  \label{QQ}
Suppose that at the $k$-th iteration of Algorithm \ref{al.pqna},  after applying the subproblem solver satisfying \eqref{Q.linear} for $r(k)$ iterations, starting with
$u_0=x_k$, we obtain solution $x_{k+1}=u_{r(k)}$. 
Let $r(k)$ satisfy 
\bequation \label{f}
{r(k)} \geq  k \log_{1/\alpha} (1/\rho^{\delta}),
\eequation
for some $\delta >1$, and $\rho$  defined in Theorem \ref{sec.linear}. 
Then 
\bequationn
Q_{H_k}( x_{k+1},x_k) -Q_k^* \leq C \rho^{k \cdot \delta},
\eequationn
holds for all $k$, with $C$ being the uniform bound on  $Q_{H_k}(x_k,x_k) - Q_k^*$, and the linear convergence of Algorithm \ref{al.pqna} is achieved. 

\etheorem
\bproof
First, assume that at the $k-$th iteration we have applied the subproblem solver for ${r(k)}$ iterations to minimize strongly convex function $Q_{H_k}$. Now, by combining $Q_{H_k}(u_0,u_0) - Q_k^* \leq C$ and \eqref{Q.linear}, we can conclude the following upper bound, so that
\bequationn  \label{alpha.C}
Q_{H_k}(u_{r(k)},u_0) - Q_k^* \leq \alpha^{r(k)} C.
\eequationn

Now, if $ \alpha^{r(k)} C \leq \epsilon_k$, we can guarantee that $u_{r(k)}$ is an $\epsilon_k$-solution of the $k$-th subproblem, so that $Q_{H_k}(u_{r(k)},u_0) \leq Q_k^* + \epsilon_k$. Now, assuming that $\rho$ is known, we can set the error rate of the $k$-th iteration as $\epsilon_k \leq C\rho^{k \cdot \delta}$, for a fixed $\delta>1$. In this case, the number of inner iterations which guarantees the $\epsilon_k$-minimizer will be
\bequationn \label{numsteps}
{r(k)} \geq k \log_{1/\alpha} (1/\rho^{\delta}).
\eequationn
\eproof

\bremark  \label{re.r}
Since subproblems are strongly convex, the required linear convergence rate for the subproblem solver, stated in \eqref{Q.linear},   can be guaranteed via some basic first-order algorithms or their accelerated variants. However, one difficulty in obtaining lower bound \eqref{f} is that it  depends on the prior knowledge of $\rho$ and $\alpha$. Consider the following simple modification of Theorem \ref{QQ}; instead of condition \eqref{f}, consider $r(k)$ satisfying 
\bequation \label{sub.r}
{r(k)} \geq k \log_{1/\alpha^\prime} ({k}/{\ell}),
\eequation
 for any given $\ell >0$ and $\alpha^\prime \in (0,1)$. Then $\epsilon_k \leq C\( {\ell}/{k} \)^k$ implies $\epsilon_k\leq C\rho^{k \cdot \delta}$, for sufficiently large $k$. 
 \eremark

In the next subsection, we extend our analysis to the case of solving subproblems via the randomized coordinate descent, where at each iteration the desired error bound related to $\epsilon_k$ is only satisfied in expectation.

\subsection{Solving Subproblems via Randomized Coordinate Descent}
As we mentioned before, in order to achieve linear convergence rate of the inexact {PQNA}, any simple first-order method (such as {PGA}) can be applied. However, as discussed in \cite{tang}, in the case when $g(x) = \lambda \|x\|_1$ and $H_k$ is sum of a diagonal and a low rank matrix, as in the case  of L-BFGS approximations, the \textit{coordinate descent} method is the most efficient approach to solve the strongly convex quadratic subproblems. 
In  this case, each iteration of coordinate descent has complexity of $\mathcal{O}(m)$, where $m$ is the memory size of L-BFGS, which is usually chosen to be less than $20$, while each iteration of a proximal gradient method has complexity of $\mathcal{O}(nm)$ and each iteration of the Newton type proximal method has complexity of $\mathcal{O}(nm^2)$. While more iterations of coordinate descent may be required to achieve the same accuracy, it tends to be the most efficient approach.
To extend our theory of the previous section and to establish the bound on the number of coordinate descent steps needed to solve each subproblem, we utilize convergence results for the \textit{randomized coordinate descent} \cite{martin}, as is done in  \cite{tang}. 

 Algorithm \ref{al.rcd} shows the framework of the randomized coordinate descent method, which can be used as a subproblem solver of Algorithm \ref{al.pqna} and is identical to the method used in \cite{tang}. In Algorithm \ref{al.rcd}, function $Q_H$ is iteratively minimized over a randomly chosen coordinate, while the other coordinates remain fixed.


\balgorithm[ht]
  \small
  \caption{Randomized Coordinate Descent Algorithm}
 \label{al.rcd}
  \balgorithmic[1]
    \smallskip
    \AlgBreak
      \State Initialize point $v \in \mathbb{R}^n$ and required number of iterations $r >0$.
      \State  Set ${u_0} \gets v$.
      \For{$l=1,2,\cdots, r-1$}
      \State Choose $j$ uniformly from $ \{1,2,\cdots,n\}$.
      \State Compute $z^* := \arg \min_{z \in \mathbb{R}^n} Q_H(u_l+ze_j,v)$.
      \State Set $u_{l+1} \gets u_l+z^*e_j$.
      \EndFor
      \State Return $u_r$.
  \ealgorithmic
\ealgorithm

In what follows, we restate \textit{Theorem 6} in \cite{martin},  which establishes linear convergence rate of the randomized coordinate descent algorithm, in expectation, to solve strongly convex problems.

\begin{theorem} \label{rcd}
Suppose we apply randomized coordinate descent for $r$ iterations, to minimize the $m$-strongly convex function $Q$ with $M$-Lipschitz gradient, to obtain the random point $u_r$.

 When $u_0$ is the initial point and $Q^* := \min_{u \in \mathbb{R}^n} Q_H(u,u_0)$, for any $r$, we have
\bequation \label{rcd.lin}
\mathbb{E}\(Q_H(u_r,u_0) - Q^*\) \leq \( 1- \frac{1- \phi_{m,M}}{n} \)^r  \(Q_H(u_0,u_0) - Q^*\),
\eequation
where $\phi_{m}$ is defined as 
\bequation \label{phi.def}
\phi_{m,M} = \begin{cases} 1-{m}/{4M} & \text {if} ~~ m \leq 2M, \\ {M}/{m} & \text {otherwise.} \end{cases}
\eequation

\end{theorem}
\bproof
The proof can be found in \cite{martin}.
\eproof

Now, we want to analyze how we can utilize the result of Theorem \ref{rcd} to achieve the linear convergence rate of inexact {PQNA}, in expectation. Toward this end, first we need the following theorem as the probabilistic extension of Theorem \ref{sec.linear}.

\btheorem \label{sec.linear.prob}
Suppose that Assumptions \ref{assum.main} and \ref{assum.strong} hold. At each iteration $k$ of the inexact {PQNA}, stated in Algorithm \ref{al.pqna}, assume that the error
$\epsilon_k$ is a nonnegaitve random variable defined on some probability space with an arbitrary distribution. Then, we have
\bequation \label{i.rate.strong}
\mathbb{E}\(F(x_k) - F(x_*)\) \leq \rho^k \( F(x_0) - F(x_*) +B_k \),
\eequation
when $\rho = 1 - {(\eta \gamma)}/{(\gamma + M)}$, and 
\bequationn
B_k := \eta \sum_{i=1}^{k} \(\mathbb{E}({\epsilon_i}) / {\rho^i} \).
\eequationn
\etheorem
\bproof
The proof is a trivial modification of  that of Theorem \ref{sec.linear}. 
\eproof

In what follows, we describe how the randomized coordinate descent method ensures the required accuracy of subproblems and consequently guarantees linear convergence of the inexact {PQNA}.

 \btheorem
 Suppose that at the $k$-th iteration of Algorithm \ref{al.pqna}, after  applying Algorithm \ref{al.rcd} 
 for  $r(k)$ iterations, starting with
$u_0=x_k$, we obtain solution $x_{k+1}=u_{r(k)}$. 
If 
 \bequationn
{r(k)} \geq k \log_{1/\alpha_n} ({k}/{\ell}), 
\eequationn
 where $\ell$ is any positive constant, $ \alpha_n = \( 1- \frac{1- \phi_{m,M}}{n} \)$ with $\phi_{m,M}$ defined in \eqref{phi.def}, and $C$ is the uniform bound on
 $Q_{H_k}(x_k,x_k) - Q_k^*$, then 
 Algorithm \ref{al.pqna}, converges linearly with constant $\rho = 1 - {(\eta \gamma)}/{(\gamma + M)}$, in expectation. 
\etheorem

\bproof
Suppose that at the $k$-th iteration of Algorithm \ref{al.pqna}, we apply $r(k)$ steps of Algorithm \ref{al.rcd} to minimize the strongly convex function $Q_{H_k}$. If $u_{r(k)}$ denotes the resulting random point, when $u_0$ is the initial point and $Q_k^* := \min_{u \in \mathbb{R}^n} Q_{H_k}(u,u_0)$, then based on Theorem \ref{rcd} we have
\bequation \label{rcd.lin2}
\mathbb{E}\(Q_{H_k}(u_{r(k)},u_0) - Q_k^*\) \leq \alpha_n^{r(k)} C,
\eequation
where $\alpha_n = \( 1- \frac{1- \phi_{m,M}}{n} \)$, with $\phi_{m,M}$ defined in \eqref{phi.def}, $Q_{H_k}(u_0,u_0) - Q_k^*$ is bounded from above by $C$. Now, based on the result of Theorem \ref{sec.linear.prob}, if $\mathbb{E}(\epsilon_k) \leq C(\ell/k)^k$ for some given positive constant $\ell$, then for sufficiently large $k$, we can guarantee that $B_k$ is uniformly bounded for all $k$, and consequently the linear convergence rate of Algorithm \ref{al.pqna}, in expectation is established. Now, by using \eqref{rcd.lin2},  $E(\epsilon_k) \leq C(\ell/k)^k$ simply follows from 
 \bequationn
{r(k)} \geq k \log_{1/\alpha_n} ({k}/{\ell}). 
\eequationn
\eproof

\bremark
The bound on the number of steps $r(k)\geq k \log_{1/\alpha_n} ({k}/{\ell})$ for randomized coordinate descent differs from the bound 
$r(k)\geq k \log_{1/\alpha} ({k}/{\ell})$ on the number of steps of a deterministic linear convergence method, such as {PGA} by the difference in constants $\alpha$ and $\alpha_n$. It can be easily shown that in the worst case $\alpha_n\approx \alpha/n$, and hence, the number of coordinate descent steps is around $n$ times larger than that of a
proximal gradient method. On the other hand, each coordinate descent step is $n$ times less expensive and in many practical cases a modest number of iterations of randomized coordinate descent is sufficient. Discussions on this can be found in \cite{martin} and \cite{tang} as well as in Section \ref{sec.result}. 
\eremark

\section{Accelerated Proximal Quasi-Newton Algorithm}\label{ch_accelerated}
We now turn  to an accelerated  variant of {PQNA}.  As we described in the introduction section, the algorithm proposed in \cite{toh} is a proximal quasi-Newton variant of FISTA, described in Algorithm  \ref{p.al1}.  In \cite{toh}, the convergence  rate of $\mathcal{O}({1}/{k^2})$ is shown under the condition that  the Hessian estimates satisfy  $0 \prec H_k \preceq H_{k-1}$, at each iteration. On the other hand, the sequence $\{ H_k \}$ is chosen so that the quadratic approximation of $f$ is an over approximation. This leads to an unrealistic setting where two possible contradictory conditions need to be satisfied and as mentioned earlier, this condition contradicts the assumptions of the original {APGA}, stated in Algorithm \ref{p.al1}.
We propose a  more general version, henceforth referred to as {APQNA},   which allows a more general sequence of $H_k$
and is based on the relaxed version of FISTA, proposed in \cite{goldfarb}, which does not impose monotonicity of the step-size parameters. Moreover, our algorithm allows more general Hessian estimates as we  explain below. 

\subsection{Algorithm Description}\label{sec.algorithm}    
The main framework of {APQNA} as stated in Algorithm \ref{p.al2} is similar to that of Algorithm \ref{p.al1}, where the simple composite quadratic approximation $Q_\mu$ was replaced by the scaled version $Q_H$, as is done in  Algorithm \ref{al.pqna}, using (partial) Hessian information.
As in the case of Algorithm \ref{al.pqna}, we assume that the approximate Hessian $H_k$ is a positive definite matrix such that $m I \preceq  H_k \preceq M I$, for some positive constants $m$ and $M$. As discussed in Remark \ref{G.bound}, it is simple to show that this condition can be satisfied for any positive $m$ and for any large enough $M$. Here, however, we will need additional much stronger assumptions on the sequence $\{H_k\}$. 
The algorithm, thus, has some additional steps compared to Algorithm \ref{p.al1} and the standard FISTA type proximal quasi-Newton algorithm proposed in \cite{toh}. Below, we present Algorithm \ref{p.al2} and  discuss the steps of each iteration in detail. 


\balgorithm[ht]
  \small
  \caption{Accelerated Proximal Quasi-Newton Algorithm}
\label{p.al2}
  \balgorithmic[1]
    \smallskip
    \AlgBreak
      \State Initialize $t_1=1,\theta_0=1, \sigma_1^0 >0$, $y_1=x_{-1}=x_0\in \mathbb{R}^n$, and positive definite matrix $H_0 \in \mathbb{R}^{ n \times n}$, and choose $ \beta \in (0,1)$.
      \For{$k=1,2,\cdots$}
      \State Define $\sigma_k:=\sigma_k^0$.
      \State Compute $p_{H_k}(y_k):=\arg \min_{u \in \mathbb{R}^n} Q_{H_k}(u,y_k)$.
      \While{$F\(p_{H_k}(y_k)\) > Q_{H_k}\(p_{H_k}(y_k),y_k\)$}
      \State Set $H_k \gets \frac{1}{\beta} H_k$.
      \State Modify $\sigma_k$ so that $\sigma_k H_k  \preceq \sigma_{k-1} H_{k-1}$.
      \State Update $\theta_{k-1}=\sigma_{k-1} / \sigma_{k}$ and recompute $t_k$ and $y_k$ using \eqref{tal1.1}-\eqref{tal1.2}.
      \State Compute $p_{H_k}(y_k):=\arg \min_{u \in \mathbb{R}^n} Q_{H_k}(u,y_k)$.
      \EndWhile
      \State Set $x_{k} \gets p_{H_k}(y_k)$.
      \State Choose {$\sigma_{k+1}^0 >0$}  and $H_{k+1}$  so that  $\sigma_{k+1}^0 H_{k+1}  \preceq \sigma_{k} H_{k}$.
      \State Define  {$\theta_k := \sigma_k / \sigma_{k+1}^0$} and compute $t_{k+1}$ and $y_{k+1}$, so that
\bsubequations \label{tal1}
  \begin{align}
 &t_{k+1} = \frac{1}{2} \(1+\sqrt{1+4 {\theta_k}t_{k}^2}\) \label{tal1.1}\\
\text{and}~~& y_{k+1} =x_{k}+\frac{t_{k}-1}{t_{k+1}}\(x_{k}-x_{k-1}\) \label{tal1.2}. 
 \end{align}
\esubequations

       \EndFor
  \ealgorithmic
\ealgorithm

 The key requirement imposed by Algorithm \ref {p.al2} on the sequence $\{H_k\}$ is that $\sigma_{k+1} H_{k+1}  \preceq \sigma_{k} H_{k}$, while   $\theta_k := {\sigma_k}/{\sigma_{k+1}}$ is used to evaluate the accelerated parameter $t_{k+1}$ through \eqref{tal1.1}. 
  During Steps 4 and 5 of iteration $k$, initial guesses for  $\sigma_{k+1}^0$ and $H_{k+1}$ are computed and  used to define $\theta_k$, which is then used to compute  $t_{k+1}$ and $y_{k+1}$. Since the approximate Hessian $H_{k+1}$ may change during Step 2 of iteration $k+1$, $\sigma_{k+1}$  may need to change as well
  in order to satisfy condition $\sigma_{k+1} H_{k+1}  \preceq \sigma_{k} H_{k}$. In particular, we may shrink the value of $\sigma_{k+1}$ and consequently will need to recompute  $\theta_k$ and, thus,  $t_{k+1}$ and $y_{k+1}$. Therefore, the  backtracking process in Step 2 of Algorithm \ref{p.al2} involves a loop which may require repeated computations of $y_k$ and hence $\nabla f(y_k)$. 
  
  \begin{remark} 
 We do not specify how to compute $H_k$ in Algorithm \ref{p.al2}, as long as it satisfies \eqref{bound.strong} and condition $\sigma_{k+1} H_{k+1}  \preceq \sigma_{k} H_{k}$. Note that Algorithm \ref{p.al2} does not allow the use of exact Hessian information at $y_{k+1}$, i.e., $H_{k+1}=\nabla^2 f(y_{k+1})$, because it is assumed that $H_{k+1}$ is computed before $y_{k+1}$ (since $y_{k+1}$ uses the value of $\sigma_{k+1}$, whose value may have to be dependent on $H_{k+1}$). However, it is possible to use $H_{k+1}=\nabla^2 f(x_k)$ in Algorithm \ref{p.al2}. To use $H_{k+1}=\nabla^2 f(y_{k+1})$, one would need to be able to compute  $\sigma_{k+1}$ before $H_{k+1}$  and somehow ensure that condition $\sigma_{k+1} H_{k+1}  \preceq \sigma_{k} H_{k}$ is satisfied. This condition can eventually be satisfied by applying similar technique to Step 2, but in that case $H_{k+1}$ will not be equal to the Hessian, but to some multiple of the Hessian, i.e.,  $\frac{1}{\beta^i}\nabla^2 f(y_{k+1})$, for some $i$. 
 
In our numerical results, we construct $H_k$ via L-BFGS and ignore condition $\sigma_{k+1} H_{k+1}  \preceq \sigma_{k} H_{k}$, since enforcing it in this case causes a very rapid decrease in $\sigma$. It is unclear, however, if a practical version of Algorithm \ref{p.al2}, based on L-BFGS Hessian approximation can be derived, which may explain why the accelerated version of our algorithm does not represent any significant advantage.
\end{remark}

One trivial choice of the matrix sequence is $H_k = \frac{1}{\mu_k} I$. In this case, the sequence of scalars $\sigma_k  =  \mu_k $, satisfies
  $\sigma_{k+1} H_{k+1}  \preceq \sigma_{k} H_{k}$, for all $k$. This choice of Hessian reduces Algorithm \ref{p.al2} to the version of {APGA} with full backtracking of the step-size parameters, proposed in \cite{goldfarb}, hence Algorithm \ref{p.al2} is the generalization of that algorithm.  Another choice for the matrix sequence is $H_k = \frac{1}{\sigma_k} H$, where the matrix $H$ is any fixed positive definite matrix. This setting of $H_k$ automatically satisfies condition $\sigma_{k+1} H_{k+1}  \preceq \sigma_{k} H_{k}$, and  Algorithm \ref{p.al2} reduces to the simplified version stated below in Algorithm \ref{p.al3}. 


\balgorithm[ht]
  \small
  \caption{Accelerated Proximal Quasi-Newton Algorithm with \\ Fixed~Hessian}
\label{p.al3}
  \balgorithmic[1]
    \smallskip
    \AlgBreak
      \State Initialize $t_1=1,\theta_0=1, \sigma_1^0 >0$, and $y_1=x_{-1}=x_0\in \mathbb{R}^n$, and choose positive definite matrix $H \in \mathbb{R}^{ n \times n}$, and $ \beta \in (0,1)$.
      \For{$k=1,2,\cdots$}
      \State Define $\sigma_k:=\sigma_k^0$.
      \State Compute $H_k = (1/\sigma_k) H$ and $p_{H_k}(y_k):=\arg \min_{u\in \mathbb{R}^n} Q_{H_k}(u,y_k)$.
      \While{$F\(p_{H_k}(y_k)\) > Q_{H_k}\(p_{H_k}(y_k),y_k\)$}
      \State Set $\sigma_k \gets {\beta} \sigma_k$.
      \State Update $\theta_{k-1}$ and recompute $t_k$ and $y_k$ using \eqref{tal2.1}-\eqref{tal2.2}.
      \State Update $H_k = (1/\sigma_k) H$.
      \State Compute $p_{H_k}(y_k):=\arg \min_{u\in \mathbb{R}^n} Q_{H_k}(u,y_k)$.
      \EndWhile
      \State Set $x_{k} \gets p_{H_k}(y_k)$.
      \State Choose {$\sigma_{k+1}^0 >0$}, define {$\theta_k := \sigma_k / \sigma_{k+1}^0$}, and compute $t_{k+1}$ and $y_{k+1}$, so that 
\bsubequations \label{tal1}
  \begin{align}
& t_{k+1} =\frac{1}{2}\(1+\sqrt{1+4 {\theta_k}t_{k}^2}\) \label{tal2.1}\\
 \text{and}~~&y_{k+1} =x_{k}+\frac{t_{k}-1}{t_{k+1}}\(x_{k}-x_{k-1}\) \label{tal2.2}. 
 \end{align}
\esubequations

      \EndFor
  \ealgorithmic
\ealgorithm

Note that, by the same logic that was used in Remark \ref{G.bound}, the number of backtracking steps at each iteration of Algorithm \ref{p.al3} is uniformly bounded. Thus, as long as the fixed approximate Hessian $H$ is positive definite, a Hessian estimate $H_k = \frac{1}{\sigma_k} H$ has positive eigenvalues bounded from above and below. In our implementation, we compute a fixed matrix $H$ by applying L-BFGS for a fixed number of iterations and then apply Algorithm \ref{p.al3}.

In the next section, we analyze the convergence properties of Algorithm \ref{p.al2}, where the approximate Hessian  $H_k$ is produced by some generic unspecified scheme.  The motivation is to be able to apply the analysis to popular and efficient Hessian approximation methods, such  as  L-BFGS. However,
in the worst case for general $H_k$, a positive lower bound for $\{\sigma_k\}$ can not be guaranteed for such a generic scheme. This observation motivates the analysis of Algorithm \ref{p.al3}, as a simplified version of Algorithm \ref{p.al2}.  It remains to be seen if some bound on $\{\sigma_k\}$ may be derived for matrices arising specifically via L-BFGS updates.

\subsection{Convergence Analysis}\label{sec.convergence}
In this section, we prove that if   the sequence $\{\sigma_k\}$ is bounded away from zero, Algorithm \ref{p.al2} achieves the same rate of convergence as {APGA}, i.e., $\mathcal{O}({1}/{k^2})$. First, we state a simple result based on  the optimality of $p_H$.

\blemma \label{e.b3}
For any $v\in \mathbb{R}^n$, there exists a subgradient of function $g$ where $\nu_g(p_H(v))\in \partial g(p_H(v))$, such that
\bequationn \label{nu}
\nabla f(v)+H(p_H(v)-v)+\nu_g(p_H(v))=0.
\eequationn
\elemma
\bproof
 The proof is followed immediately from the optimality condition of the convex optimization problem \eqref{e.b2}. 
\eproof 

Now, we can show the following lemma,  which bounds the change in the objective function $F$ and is a simple extension of \textit{Lemma 2.3} in \cite{beck}. 

\blemma \label{e.b4}
Let $v \in \mathbb{R}^n$ and $H \succ 0$ be such 
\bequation \label{e.b5}
F\(p_H(v)\) \leq Q_H\(p_H(v),v\),
\eequation 
holds for a given $v$, then for any $x \in \mathbb{R}^n$
\bequationn \label{e.b6}
F(x) - F(p_H(v)) \geq \frac{1}{2} \|p_H(v) - v \|_H^2 + \langle v-x, p_H(v)-v \rangle_H.
\eequationn
\elemma
\bproof
From \eqref{e.b5}, we have
\bequation \label{e.b7}
F(x) - F(p_H(v)) \geq F(x)-Q_H(p_H(v),v).
\eequation
Now, based on the convexity of functions $f$ and $g$, we have
\bequationn
\baligned
&f(x) \geq f(v) +  \langle \nabla f(v), x-v \rangle \\
\text{and}~~&g(x) \geq g(p_H(v)) +  \langle \nu_g(p_H(v)), x-p_H(v) \rangle,
\ealigned
\eequationn
where $\nu_g(p_H(v))$ is defined in Lemma \ref{e.b3}. Summing the above inequalities yields
\bequation \label{e.b8}
F(x) \geq f(v) +  \langle \nabla f(v), x-v \rangle +g(p_H(v)) +  \langle \nu_g(p_H(v)), x-p_H(v) \rangle.
\eequation
Using \eqref{e.b1} and \eqref{e.b8} in \eqref{e.b7} yields
\bequationn
\baligned
F(x) - F(p_H(v)) \geq &~ f(v) +  \langle \nabla f(v), x-v \rangle +g(p_H(v)) + \langle \nu_g(p_H(v)), x-p_H(v) \rangle \\
&-f(v) - \langle \nabla f(v), p_H(v)-v \rangle - \frac{1}{2} \|p_H(v)-v\|_H^2 - g(p_H(v)) \\
 =& - \frac{1}{2} \|p_H(v)-v\|_H^2 + \langle x-p_H(v), \nabla f(v)+\nu_g(p_H(v))\rangle\\
=& - \frac{1}{2} \|p_H(v)-v\|_H^2 + \langle x-p_H(v), H(v-p_H(v)) \rangle \\
=& - \frac{1}{2} \|p_H(v)-v\|_H^2 + \langle x-p_H(v), H(v-p_H(v)) \rangle \\
& + \langle v-p_H(v), v-p_H(v)\rangle_H - \langle v-p_H(v), v-p_H(v)\rangle_H \\
=&~ \frac{1}{2} \|p_H(v)-v\|_H^2 + \langle v-x, p_H(v) - v\rangle_H.
\ealigned
\eequationn
\eproof 

The following result is a simple corollary of Lemma \ref{e.b4}. 

\bcorollary
Let $v \in \mathbb{R}^n$ and $H \succ 0$ be such that
\bequationn 
F(p_H(v)) \leq Q_H(p_H(v),v),
\eequationn 
 then for any $x \in \mathbb{R}^n$
\bequation \label{new.main}
\begin{aligned}
2(F(x)-F(p_H(v))) &\geq  \|p_H(v)-v \|_H^2+ 2 \langle p_H(v)-v, v-x \rangle_H,\\
&= \| p_H(v)-x \|^2_H- \|v-x\|^2_H. 
\end{aligned}
\eequation
\ecorollary
\bproof
The result immediately follows by applying the following identity 
\bequation \label{py}
\|b-a\|^2+2(b-a)^T(a-c)=\|b-c\|^2-\|a-c\|^2.
\eequation
to  Lemma \ref{e.b4} with 
\bequationn
a:=H^{\frac{1}{2}}v,~~~ b:=H^{\frac{1}{2}} p_H(v),~~~c:=H^{\frac{1}{2}}x. 
\eequationn
\eproof

The next lemma states the key properties which are used in the convergence analysis.
\blemma  \label{lem:basicprop}
At each iteration of  Algorithm \ref{p.al2}, the following relations hold
\bsubequations \label{assum2}
  \begin{align}
 & \sigma_k H_k  \succeq \sigma_{k+1} H_{k+1} \label{as2} \\
 \text{and}~~ &\sigma_k t_k^2  \geq \sigma_{k+1} t_{k+1} (t_{k+1}-1).  \label{as1}
  \end{align}
\esubequations
\elemma
\bproof
The proof follows trivially from the conditions in Algorithm \ref{p.al2} and the fact that $\theta_k \leq {\sigma_k}/{\sigma_{k+1}}$.
\eproof
 Now, using this lemma and  previous results we derive the key property of the iterations of {APQNA}. 

\blemma \label{lemma3}
For all $k\geq 1$  for Algorithm \ref{p.al2}  we have
 \bequationn \label{3}
 2\sigma_kt_k^2v_k+\sigma_ku_k^TH_ku_k \geq 2\sigma_{k+1}t_{k+1}^2v_{k+1}+\sigma_{k+1}u_{k+1}^TH_{k+1}u_{k+1},
 \eequationn
 where $v_k=F(x_k)-F(x_*)$ and $u_k=t_kx_k - (t_k-1)x_{k-1}-x_*$.
\elemma

\bproof
In \eqref{new.main}, by setting $v=y_{k+1}$, $p_H(v)=x_{k+1}$, $H=H_{k+1}$, and $x=x_k$ and then by multiplying the resulting inequality by $\sigma_{k+1}(t_{k+1}-1)$, we will have
\bequationn \label{hi1}
\begin{aligned}
& 2\sigma_{k+1}(t_{k+1}-1)(v_k-v_{k+1}) \\
\geq & ~ (t_{k+1}-1)(x_{k+1}-y_{k+1})^T\sigma_{k+1}H_{k+1}(x_{k+1}-y_{k+1})\\
&+2(t_{k+1}-1)(x_{k+1}-y_{k+1})^T\sigma_{k+1}H_{k+1}(y_{k+1}-x_k).
\end{aligned}
\eequationn
On the other hand, in \eqref{new.main}, by setting $x=x_*$ and multiplying it by $\sigma_{k+1}$, we have
\bequationn
\begin{aligned}
-2\sigma_{k+1}v_{k+1} \geq&~ (x_{k+1}-y_{k+1})^T\sigma_{k+1}H_{k+1}(x_{k+1}-y_{k+1})\\
&+2(x_{k+1}-y_{k+1})^T\sigma_{k+1}H_{k+1}(y_{k+1}-x_*).
\end{aligned}
\eequationn
By adding these two inequalities, we have
\bequationn
\begin{aligned}
& 2\sigma_{k+1}((t_{k+1}-1)v_k-t_{k+1}v_{k+1}) \\
 \geq &~t_{k+1} (x_{k+1}-y_{k+1})^T \sigma_{k+1}H_{k+1}(x_{k+1}-y_{k+1})\\
&+2(x_{k+1}-y_{k+1})^T \sigma_{k+1}H_{k+1}(t_{k+1}y_{k+1}-(t_{k+1}-1)x_k-x_*).
\end{aligned}
\eequationn
Multiplying the last inequality by $t_{k+1}$ and applying inequality \eqref{as1} give
\bequationn
\begin{aligned}
& 2(\sigma_k t_k^2v_k-\sigma_{k+1}t_{k+1}^2v_{k+1}) \\
\geq& ~ t_{k+1}^2 (x_{k+1}-y_{k+1})^T \sigma_{k+1}H_{k+1}(x_{k+1}-y_{k+1})\\
&+2t_{k+1}(x_{k+1}-y_{k+1})^T \sigma_{k+1}H_{k+1}(t_{k+1}y_{k+1}-(t_{k+1}-1)x_k-x_*).
\end{aligned}
\eequationn
 By applying \eqref{py} with
\bequationn
\begin{aligned}
&a:=\sqrt{\sigma_{k+1}} H_{k+1}^{\frac{1}{2}}t_{k+1}y_{k+1},~~b:=\sqrt{\sigma_{k+1}} H_{k+1}^{\frac{1}{2}}t_{k+1}x_{k+1},\\
&c:=\sqrt{\sigma_{k+1}} H_{k+1}^{\frac{1}{2}}((t_{k+1}-1)x_k+x_*),
\end{aligned}
\eequationn
the last inequality can be written as
\bequationn \label{hi3}
\begin{aligned}
&2(\sigma_k t_k^2v_k-\sigma_{k+1}t_{k+1}^2v_{k+1}) \\
& \geq ~ \|\sqrt{\sigma_{k+1}}H_{k+1}^{\frac{1}{2}}t_{k+1}x_{k+1}-\sqrt{\sigma_{k+1}} H_{k+1}^{\frac{1}{2}}((t_{k+1}-1)x_k+x_*)\|^2\\
&-\|\sqrt{\sigma_{k+1}} H_{k+1}^{\frac{1}{2}}t_{k+1}y_{k+1}-\sqrt{\sigma_{k+1}} H_{k+1}^{\frac{1}{2}}((t_{k+1}-1)x_k+x_*)\|^2.
\end{aligned}
\eequationn
Hence, by using the definition of $y_{k+1}$ and $u_k$, we have
\bequationn \label{res}
\begin{aligned}
2(\sigma_k t_k^2v_k-\sigma_{k+1}t_{k+1}^2v_{k+1}) &\geq u_{k+1}^T\sigma_{k+1}H_{k+1}u_{k+1}-u_k^T\sigma_{k+1}H_{k+1}u_k.
\end{aligned}
\eequationn
Now, based on \eqref{as2}, we have
\bequationn
 u_k^T\sigma_{k}H_{k}u_k \geq u_k^T\sigma_{k+1}H_{k+1}u_k,
\eequationn
which implies
\bequationn
2(\sigma_k t_k^2v_k-\sigma_{k+1}t_{k+1}^2v_{k+1}) \geq u_{k+1}^T\sigma_{k+1}H_{k+1}u_{k+1}-u_k^T\sigma_{k}H_{k}u_k.
\eequationn
\eproof 
Now, we are ready to state and prove the  convergence  rate result. 

\btheorem \label{theorem1}
The sequence of iterates $x_k$, generated by Algorithm \ref{p.al2}, satisfies
\bequationn \label{6}
F(x_k)-F(x_*) \leq \frac{\|x_0-x_*\|^2}{2 \sigma_k t_k^2}.
\eequationn
\etheorem

\bproof By setting $t_1=1$, using the definition of $u_k$ at $k=1$, which is $u_1=x_1-x_*$, and also considering the positive definiteness of $H_k$ for all $k\geq 1$, it follows from Lemma \ref{lemma3} that
\bequation \label{7}
2\sigma_k t_k^2v_k \leq 2\sigma_k t_k^2v_k +\sigma_k u_k^TH_ku_k\leq 2\sigma_1t_1^2v_1+(x_1-x_*)^T \sigma_1H_1 (x_1-x_*).
\eequation
Setting $x=x_*$, $v=y_1=x_0$, $p_H(v)=x_1$, $t_1=1$, and $H=H_1$ in \eqref{new.main} implies
\bequationn
-2v_1 \geq (x_1-x_*)^T H_1 (x_1-x_*)-(x_0-x_*)^T H_1 (x_0-x_*).
\eequationn
Multiplying the above  by $\sigma_1$ gives
\bequationn
2\sigma_1v_1 + (x_1-x_*)^T \sigma_1H_1 (x_1-x_*)\leq(x_0-x_*)^T \sigma_1H_1 (x_0-x_*).
\eequationn
By using inequality \eqref{7}, we  have
\bequationn
2\sigma_k t_k^2v_k \leq (x_0-x_*)^T \sigma_1H_1 (x_0-x_*).
\eequationn
Finally, by setting $\sigma_1=1$ and $H_1=I$, we obtain
\bequationn
v_k \leq \frac{\|x_0-x_*\|^2}{2\sigma_k t_k^2},
\eequationn
which completes the proof.
\eproof 

Now, based on the result of Theorem \ref{theorem1}, in order to obtain the rate of convergence of $\mathcal{O}({1}/{k^2})$ for Algorithm \ref{p.al2}, it is sufficient to show that 
\bequationn \label{9}
 \sigma_k t_k^2 \geq \psi k^2,
\eequationn
for some constant $\psi >0$. The next result is a simple consequence of the relation  \eqref{as1}, or equivalently \eqref{tal1.1}.

\blemma \label{lemma4}
The sequence $\{\sigma_k\}$  generated by Algorithm \ref{p.al2} satisfies
\bequationn
\sigma_k t_k^2 \geq \left (\frac{\sum_{i=1}^{k} \sqrt{\sigma_i}}{2}\right )^2.
\eequationn
\elemma

\bproof 
We can prove this lemma by using induction. Trivially, for $k=1$, since $t_1=1$, the inequality holds. As the induction assumption, assume that for $k >1$, we have $\sigma_k t_k^2 \geq \left (\frac{\sum_{i=1}^{k} \sqrt{\sigma_i}}{2}\right)^2$. Since \eqref{tal1.1} holds for all $k$, it follows that
\bequationn
\begin{aligned}
t_{k+1}=\frac{1}{2}+\sqrt{\frac{1}{4}+(\frac{\sigma_{k}}{\sigma_{k+1}})t_{k}^2}\geq\frac{1}{2}+\sqrt{\frac{\sigma_{k}}{\sigma_{k+1}}} t_{k}.
\end{aligned}
\eequationn
Multiplying by $\sqrt{\sigma_{k+1}}$ implies
\bequationn
\sqrt{\sigma_{k+1}} t_{k+1} \geq \frac{\sqrt{\sigma_{k+1}}}{2}+\sqrt{\sigma_k}t_k.
\eequationn
Finally, by using induction assumption, we will have have
\bequationn
\sqrt{\sigma_{k+1}} t_{k+1} \geq \frac{\sqrt{\sigma_{k+1}}}{2}+\frac{\sum_{i=1}^{k} \sqrt{\sigma_i}}{2}=\frac{\sum_{i=1}^{k+1} \sqrt{\sigma_i}}{2}. 
\eequationn 
\eproof 

Hence, if  we assume that the sequence $\{\sigma_k\}$ is bounded below by a positive constant $\underline{\sigma}$, i.e., $\sigma_k \geq \underline{\sigma}$, we can 
establish the desired bound on  $\sigma_k t_k^2$, as stated in the following theorem.

\btheorem \label{theorem3}
If for all iterations of  Algorithm \ref{p.al2} we have $\sigma_k \geq \underline{\sigma}$, then for all $k$
\bequation \label{11}
F(x_k)-F(x_*) \leq \frac{2 \|x_0-x_*\|^2}{\underline{\sigma} k^2}.
\eequation
\etheorem

\bproof
Under the assumption  $\sigma_k \geq \underline{\sigma}$, we will have
\bequationn
\(\frac{\sum_{i=1}^k \sqrt{\sigma_i}}{2}\)^2 \geq \frac{k^2  \underline{\sigma}}{4},
\eequationn
and consequently, by using Lemma \ref{lemma4}, we obtain
\bequationn
\sigma_k t_k^2 \geq \frac{k^2 \underline{\sigma}}{4}.
\eequationn
Then, by using Theorem \ref{theorem1}, we  have the desired rate of convergence of $\mathcal{O}({1}/{k^2})$ as stated in \eqref{11}.
\eproof

The assumption of the existence of a bounded sequence $\{\sigma_k\}$ such that $\sigma_k \geq \underline{\sigma}$ and \eqref{as1} holds may not be satisfied when we use a general approximate Hessian. To illustrate this, consider the following simple sequence of matrices:
\bequationn
H_{2k}=\left [ \begin{array}{cc} 10~ & 0\\ 0~ & 1\end{array} \right ]~~~\text{and}~\quad H_{2k+1}=\left [ \begin{array}{cc} 1~& 0\\ 0~ & 10\end{array} \right ]. \quad  
\eequationn
Clearly, $\sigma_{2k+1}\leq \sigma_{2k}/10$ and $\sigma_{2k}\leq \sigma_{2k-1}/10$, and hence $\sigma_k\leq 10^{-k}$. In this case, based on the result of Theorem \ref{theorem1}, we cannot guarantee any convergence result. Some convergence result can still be attained, when $\sigma_k\to 0$, for example, 
if  $\sigma_k\geq {\underline{\sigma}}/{k}$, as we show in the following relaxed version of Theorem \ref{theorem3}.
\btheorem \label{theorem3.1}
If for all iterations of  Algorithm \ref{p.al2} we have $\sigma_k \geq {\underline{\sigma}}/{k}$, then for all $k$
\bequation \label{11.1}
F(x_k)-F(x_*) \leq \frac{2 \|x_0-x_*\|^2}{\underline{\sigma} k}.
\eequation
\etheorem

\bproof
From  $\sigma_k \geq {\underline{\sigma}}/{k}$, we will have
\bequationn
\(\frac{\sum_{i=1}^k \sqrt{\sigma_i}}{2}\)^2 \geq \frac{k  \underline{\sigma}}{4},
\eequationn
and consequently, by using Lemma \ref{lemma4}, we obtain
\bequationn
\sigma_k t_k^2 \geq \frac{k \underline{\sigma}}{4}.
\eequationn
Then, by using Theorem \ref{theorem1}, we have \eqref{11.1}.
\eproof 

The above theorem shows that  if $\sigma_k$ converges to zero, but not faster than $1/k$, then our {APQNA} method may loose its accelerated rate of convergence, but still converges at least at the same rate as {PQNA}. Establishing lower bounds of $\sigma_k$ for different choices of Hessian estimates is a nontrivial task and is the subject of future research. As we will demonstrate in our computational section, {APQNA} with L-BFGS Hessian approximation does not seem to have any practical advantage over its nonaccelerated counterpart, however it is clearly convergent.

We can establish the accelerated rate of Algorithm \ref{p.al3}, since in this case we can guarantee a lower bound on $\sigma_k$, due to the restricted nature of the $H_k$ matrices.

\blemma \label{bound.sigma}
In Algorithm \ref{p.al3}, let $mI\preceq H$, then $\sigma_k\geq {\beta m}/{L}$ and hence the  convergence rate of $\mathcal{O}({1}/{k^2})$ is achieved.
\elemma
\bproof
In Algorithm \ref{p.al3}, we define $H_k = \frac{1}{\sigma_k} H$. 
The sufficient decrease condition $F(p_{H_k}(y_k)) \leq Q_{H_k}(p_{H_k}(y_k), y_k)$, is satisfied for any $H_k \succeq L I$, 
hence it is satisfied for any $H_k = \frac{1}{\sigma_k} H$ with $\sigma_k \leq { m}/{L}$.  By the mechanism of Step 3 in  Algorithm \ref{p.al3}, we observe that for all $k$, we have $\sigma_k \geq {\beta m}/{L}$.
Let us note now that Algorithm \ref{p.al3} is a special case of Algorithm   \ref{p.al2}, hence all the above results, in particular Theorem \ref{theorem1} and Lemma \ref{lemma4} hold. 
 Consequently, based on Theorem \ref{theorem3}, the desired convergence rate of $\mathcal{O}({1}/{k^2})$ for Algorithm \ref{p.al3} is obtained.
\eproof

\bremark\label{r.inexact}
We have studied only the exact variant of {APQNA} in this section. Incorporating inexact subproblem solutions, as was done for {APQNA} in the previous section, is relatively straightforward following the techniques for inexact {APGA}, \cite{mark}. It is easy to show that if the exact algorithm has the accelerated convergence rate, then the inexact counterpart, with subproblems solved by a linearly convergent method, such as randomized coordinate descent, inherits this convergence rate.  However, using the relaxed sufficient decrease condition does not apply  here as it does not preserve the accelerated convergence rate. 
\eremark

In the next section, we present the numerical results comparing the performance of Algorithm   \ref{p.al2} and Algorithm   \ref{p.al3} to their nonaccelerated counterparts, to see how much practical acceleration is achieved.  

\section{Numerical Experiments}\label{sec.result}

In this section, we investigate the practical performance of several algorithms discussed in this work, applied to the sparse logistic regression problem
\bequationn \label{logreg1}
\min_{w} \{ F(w) :=  \frac{1}{m} \sum_{i=1}^{m} \log(1+\exp(-y_i \cdot w^Tx_i))+ \lambda \|w\|_1,~ w\in \mathbb{R}^{n}\},
\eequationn
where $f(w) = \frac{1}{m}\sum_{i=1}^{m} \log(1+\exp(-y_i \cdot w^Tx_i))$ is the average logistic loss function and $g(w) = \lambda \|w\|_1$, with $\lambda > 0$, is the $\ell_1$-regularization function. 
The input data for this problem is a set of $m$ training data points,  $x_i \in \mathbb{R}^n$, and  corresponding labels $y_i \in \{-1,+1\}$, for $i=1,2,\dots,m$.

The algorithms that we compare here are as follows:
\bitemize
\item Accelerated Proximal Gradient Algorithm ({APGA}), proposed in \cite{beck}, (also known as FISTA),
\item Proximal Quasi-Newton Algorithm with Fixed Hessian approximation (call it {PQNA-FH}),
\item Accelerated Proximal Quasi-Newton Algorithm with Fixed Hessian approximation (call it {APQNA-FH}),
\item Proximal Quasi-Newton Algorithm with L-BFGS Hessian approximation (call it {PQNA-LBFGS}), proposed in \cite{tang}, and
\item Accelerated Proximal Quasi-Newton Algorithm with L-BFGS Hessian approximation (call it {APQNA-LBFGS}).
\eitemize
In {PQNA-FH} and {APQNA-FH}, we set $H_k = \frac{1}{\sigma_k} H$, where $H$ is a positive definite matrix computed via applying L-BFGS updates over the first few iterations of the algorithm which then is fixed for all remaining iterations. 
On the other hand, {PQNA-LBFGS} and {APQNA-LBFGS} employ  the L-BFGS updates to compute Hessian estimates throughout the algorithm.
In all of the above algorithms, we use the coordinate descent scheme, as described in \cite{tang}, to solve the subproblems inexactly.
According to the theory in \cite{tang}, {PQNA-FH} and {PQNA-LBFGS}  converge at the rate of $\mathcal{O}({1}/{k})$. If $f$ is strongly convex (which depends on the problem data), then according to Theorem \ref{sec.linear},{PQNA-FH} and {PQNA-LBFGS} converge at a linear rate. 
By Lemma \ref{bound.sigma},  in {APQNA-FH}, condition $\sigma_{k} H_{k}  \preceq \sigma_{k-1} H_{k-1}$ holds automatically and the algorithm converges at the rate of $\mathcal{O}({1}/{k^2})$. 
On the other hand, for {APQNA-LBFGS}, condition $\sigma_{k} H_{k}  \preceq \sigma_{k-1} H_{k-1}$  has to be enforced. 
We have tested various implementations that ensure this condition and none have produced a practical approach. We then chose to 
set $\theta_k = 1$ and relax the condition $\sigma_{k} H_{k}  \preceq \sigma_{k-1} H_{k-1}$. The resulting algorithm is practical and is empirically convergent, but as we will see does not provide an improvement over {PQNA-LBFGS}. 

Throughout all of our experiments, we initialize the algorithms with $w_0 = \bold{0}$ and we set the regularization parameter $\lambda = 10^{-3}$. Each algorithm terminates whenever $\| (\partial F(x_k))_{\text{min}}\|_{\infty} \leq 10^{-5} \| \(\partial F(x_0)\)_{\text{min}}\|_{\infty} $. In terms of the stopping criteria of subproblems solver at $i$-th iteration, we performed the coordinate descent method for $r(i)$ steps, so that $r(i) > \min(10^3, i/3)$, as long as the generated step is longer than $10^{-16}$. In {APQNA-FH} and {PQNA-FH}, in order to construct the fixed matrix $H$, we apply the L-BFGS scheme by using the information from the first $\bar k$ (with $\bar k$ chosen between $1$ and $10$) iterations and then use that fixed matrix through the rest of the algorithm. Finally, to construct the sequence $\{\sigma_k\}$, we set $\sigma_0 =1$ and $\sigma_{k+1}^0=1.015 \sigma_k$. The information on the  data sets used in our tests is summarized in Table \ref{p.table1}. These data sets are available through UCI machine learning repository \footnote{\url{http://archive.ics.uci.edu/ml/}}.

\btable[H]  
 \centering
 \captionsetup{justification=centering}
  \caption{Data information, dimension (d) and number of data points (N).}
   \label{p.table1}
  \begin{tabular}{ c|c|c|cr}
Instance &$d$&$N$&Description\\ 
\hline
\abovespace
a9a &123 & 32561& census income dataset\\ 
mnist&782&100000&handwritten digit recognition\\ 
connect-4&126&10000&win versus loss recognition\\ 
HAPT&561& 7767&human activities and postural transitions recognition\\ 
\end{tabular}
\centering
\etable
The algorithms are implemented in MATLAB R2014b and computations were performed on the COR@L computational cluster of the ISE department at  Lehigh, consisting of 16-cores AMD Operation, 2.0 GHz nodes with 32 Gb of memory.

First, in order to demonstrate  the effect of using even limited Hessian information within an  accelerated method, we compared the performance of {APQNA-FH} and {APGA}, both in terms of the number of iterations and the total solution time,  see the results in Table \ref{p.table2}.

\btable[H]     
\small
 \centering
 \captionsetup{justification=centering}
 \caption{ {APQNA-FH} vs. {APGA } in terms of function value (Fval), number of iterations (iter) and total solution time (time) in seconds.} 
  \label{p.table2}
  \begin{tabular}{ c|cc|cc|ccc}  
\multicolumn{8}{c}{\textbf{a9a}} \\ \hline  
\abovespace
 Algorithm&iter&Fval& iter& Fval&iter&Fval &time\\ \hline
\abovespace
 APGA&40&3.4891e-01&80&3.4730e-01&862&3.4703e-01&1.95e+01\\ 
\abovespace
APQNA-FH&40&3.4706e-01&80&3.4703e-01&121&3.4703e-01&5.52e+00\\ 
\multicolumn{8}{c}{\vspace{-0.1in}} \\
\multicolumn{8}{c}{\textbf{mnist}} \\ \hline
\abovespace
 Algorithm&iter&Fval& iter&Fval&iter&Fval&time\\ \hline
\abovespace
APGA&48&9.1506e-02&96&9.0206e-02&1202&8.9695e-02&5.13e+02\\ 
\abovespace
APQNA-FH&48&8.9754e-02&96&8.9699e-02&144&8.9695e-02&9.91e+01\\ 
\multicolumn{8}{c}{\vspace{-0.1in}} \\
\multicolumn{8}{c}{\textbf{connect-4}} \\ \hline
\abovespace
 Algorithm&iter&Fval& iter&Fval&iter&Fval&time\\ \hline
 \abovespace
APGA&92&3.8284e-01&184&3.7777e-01&3045&3.7682e-01&4.65e+01\\ 
\abovespace
APQNA-FH&92&3.7701e-01&184&3.7683e-01&278&3.7682e-01&2.05e+01\\ 
\multicolumn{8}{c}{\vspace{-0.1in}} \\
\multicolumn{8}{c}{\textbf{HAPT}} \\ \hline
\abovespace
 Algorithm&iter&Fval& iter& Fval&iter&Fval&time \\ \hline
 \abovespace
APGA&222&8.5415e-02&444&7.7179e-02&13293&7.1511e-02&1.38e+03\\ 
\abovespace
APQNA-FH&222&7.2208e-02&444&7.1524e-02&677&7.1511e-02&1.53e+02\\ 
\end{tabular}
\centering
\etable

Based on the results shown in Table \ref{p.table2}, we conclude that {APQNA-FH} consistently dominates the {APGA}, both in terms of the number of function evaluations and also in terms of the total solution time. 

It is worth mentioning that although in terms of computational effort, each iteration of {APGA} is  cheaper than each iteration of {APQNA-FH}, the total solution time of {APQNA-FH} is significantly less than {APGA}, due to the smaller number of iterations of {APQNA-FH} compared to {APGA}.

\begin{figure}[ht!] 
\small
\begin{subfigure}{.51\textwidth}
  \centering
  \includegraphics[width=1\linewidth]{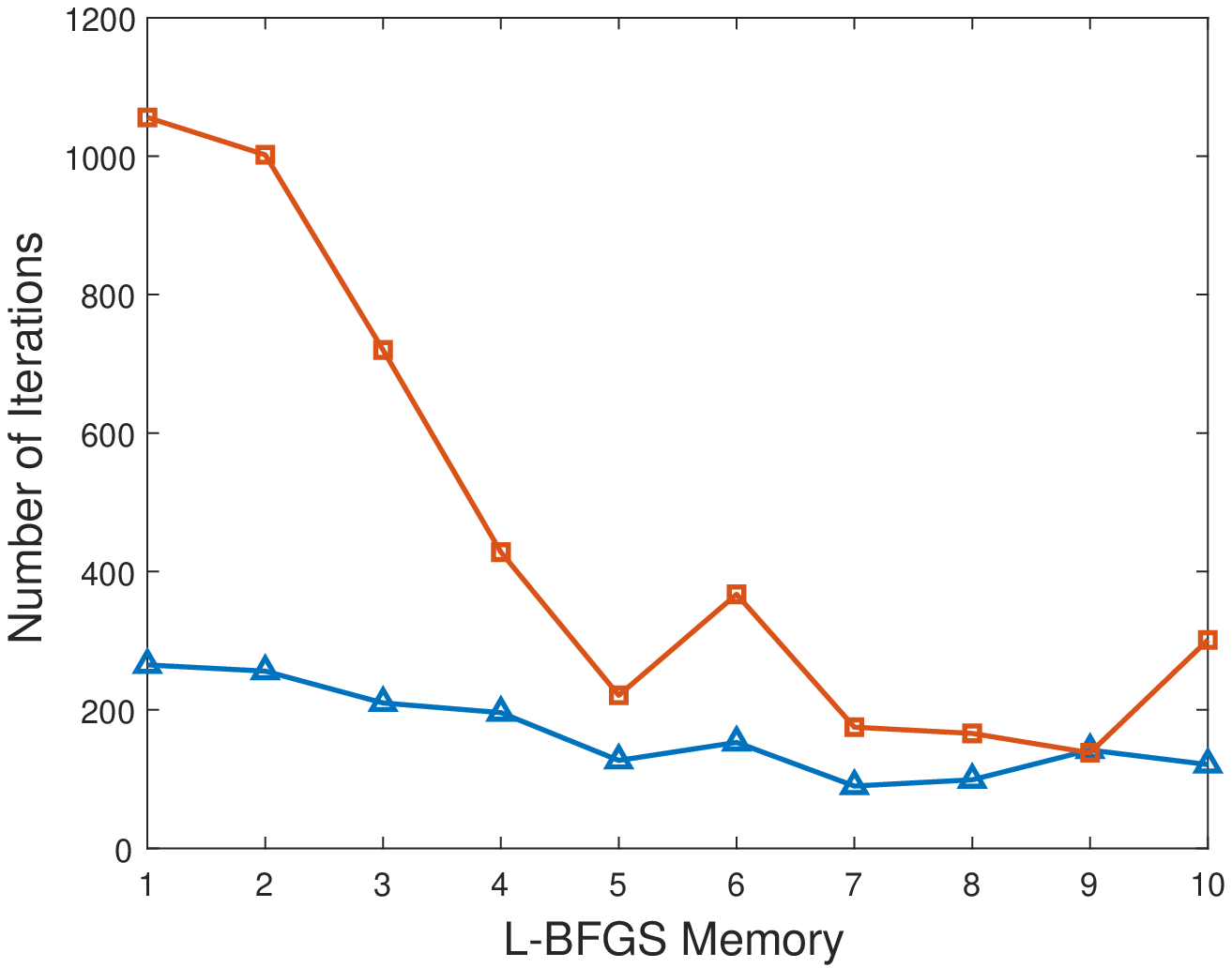}
  \caption{a9a (d=123, m=32561)}
\end{subfigure}%
\begin{subfigure}{.51\textwidth}
  \centering
  \includegraphics[width=1\linewidth]{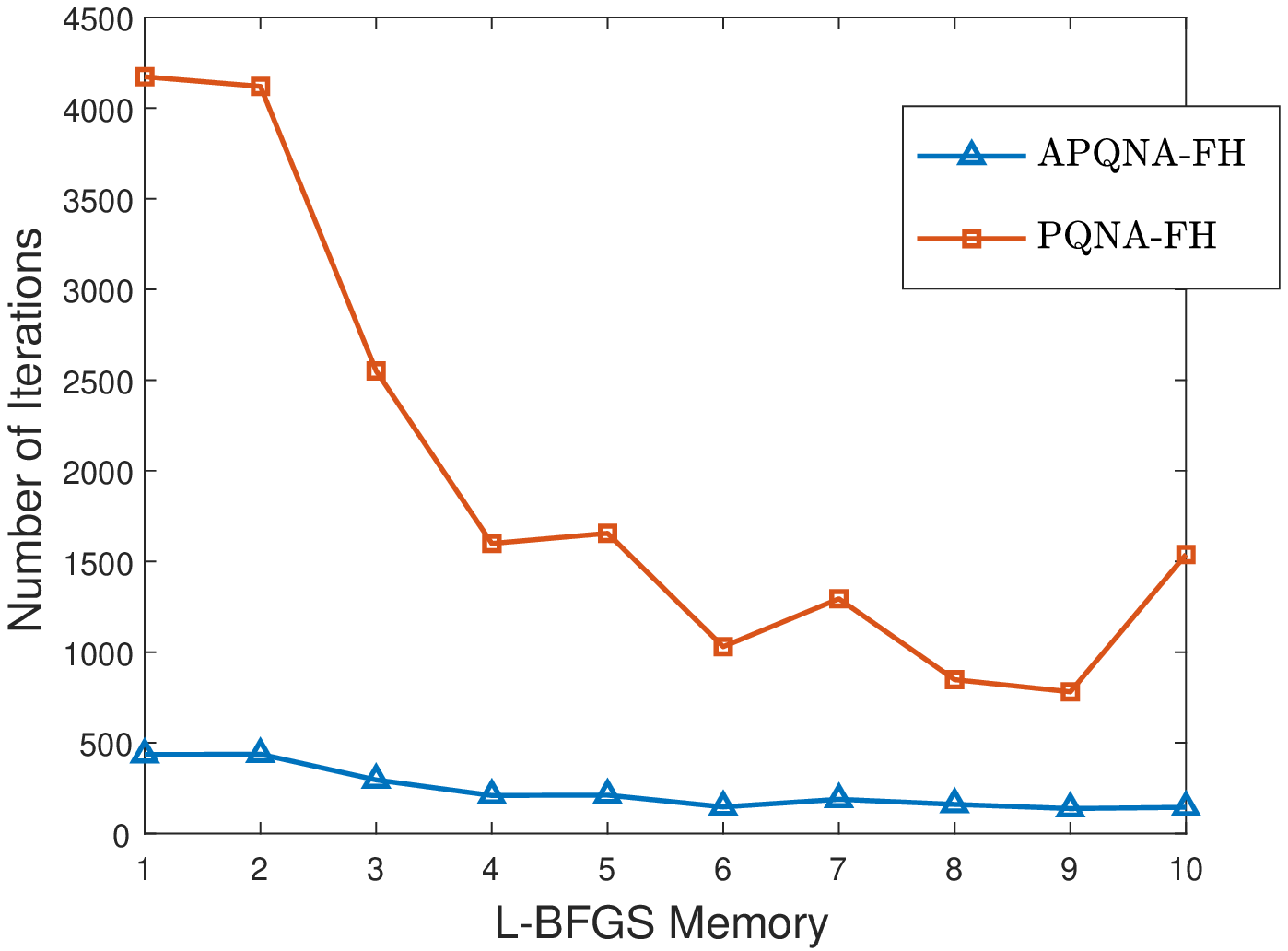}
  \caption{mnist (d=782, m=100000)}
\end{subfigure}
\begin{subfigure}{.51\textwidth}
  \centering
  \includegraphics[width=1\linewidth]{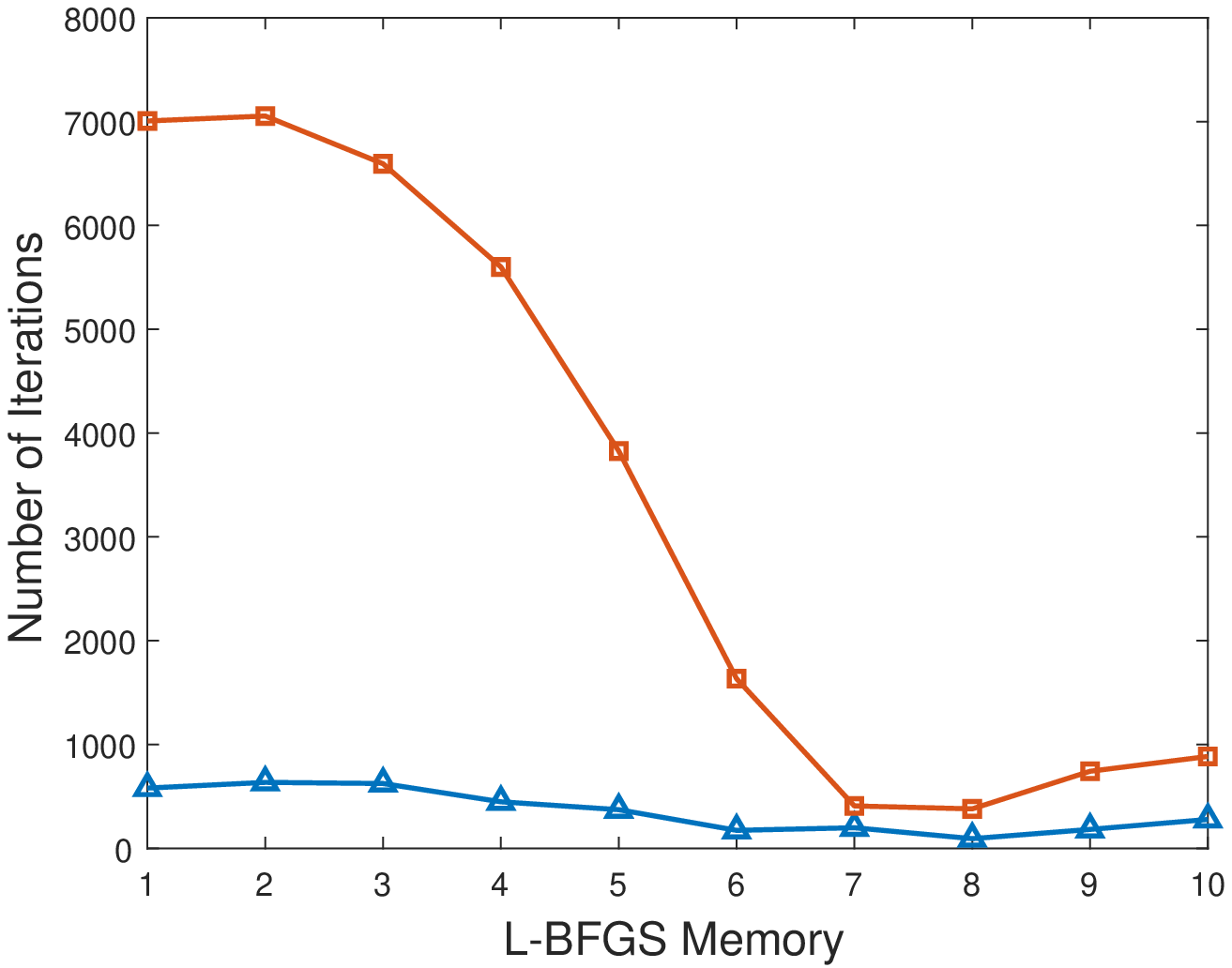}
  \caption{connect-4 (d=126, m=10000)}
\end{subfigure}
\begin{subfigure}{.51\textwidth}
  \centering
  \includegraphics[width=1\linewidth]{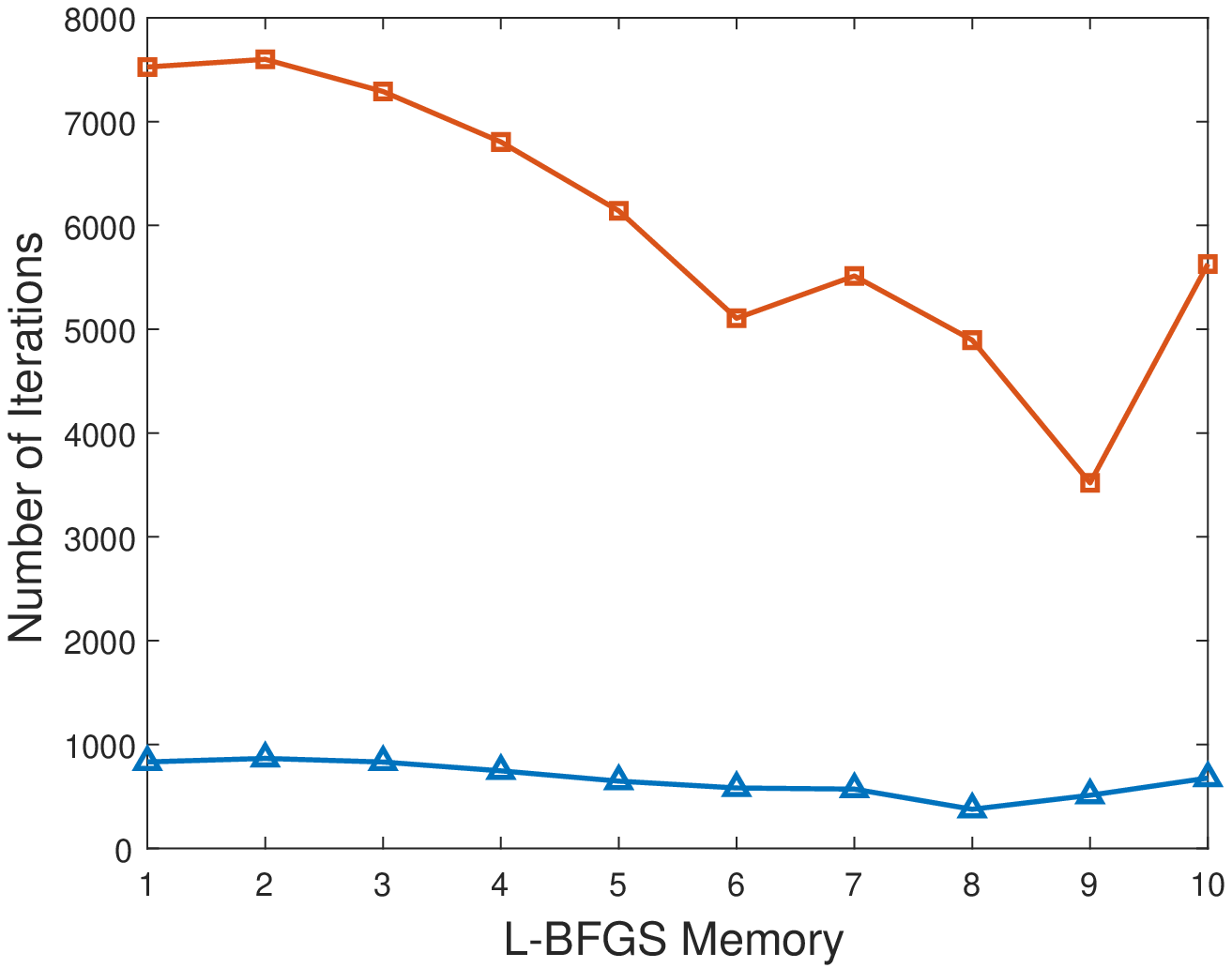}
  \caption{HAPT  (d=561, m=7767)}
\end{subfigure}
\caption{{APQNA-FH} vs. {PQNA-FH} in terms of number of iterations.}
\label{p.fig2}
\end{figure}

\begin{figure} [ht!]
\begin{subfigure}{.51\textwidth}
  \centering
  \includegraphics[width=1\linewidth]{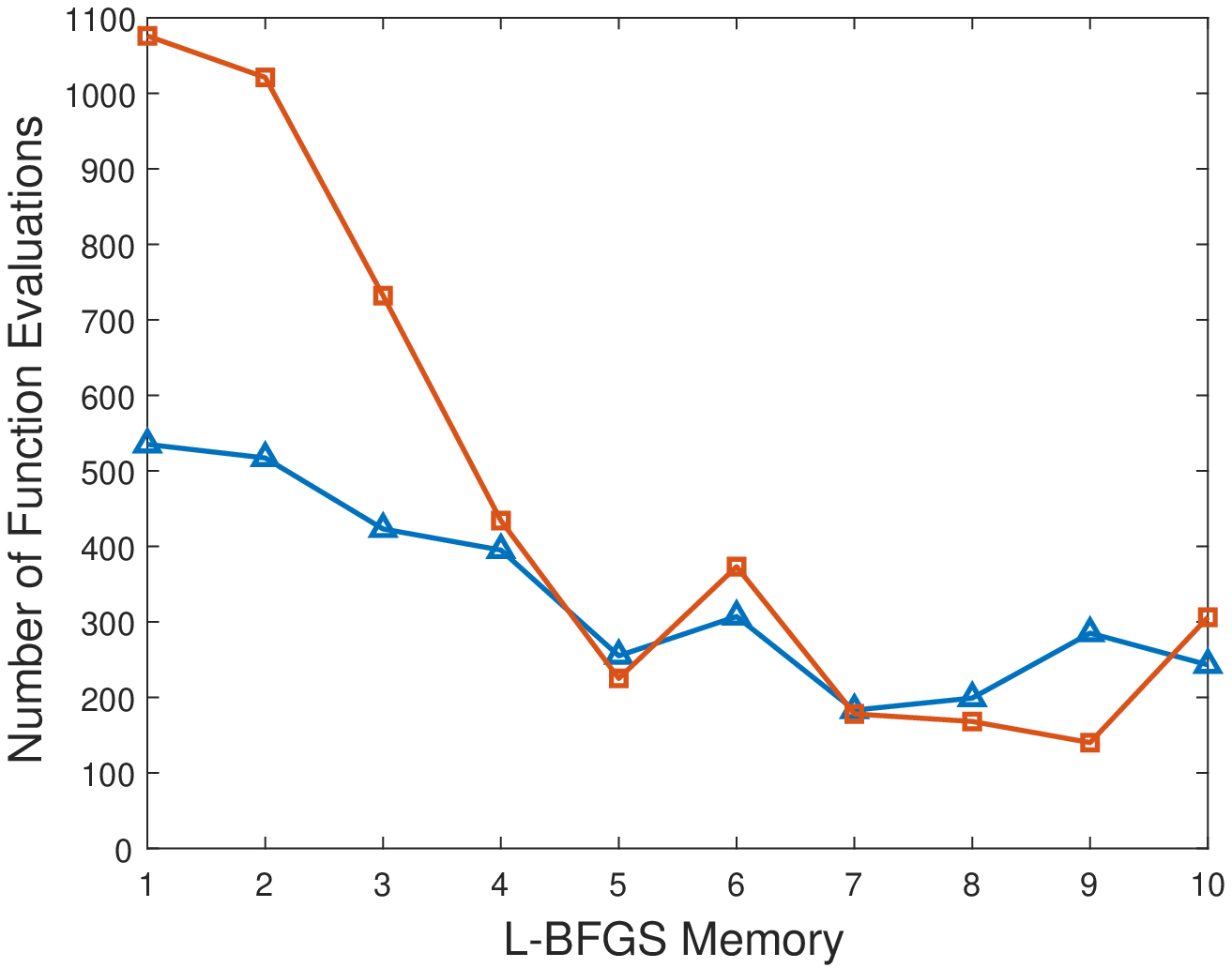}
  \caption{a9a (d=123, m=32561)}
\end{subfigure}%
\begin{subfigure}{.51\textwidth}
  \centering
  \includegraphics[width=1\linewidth]{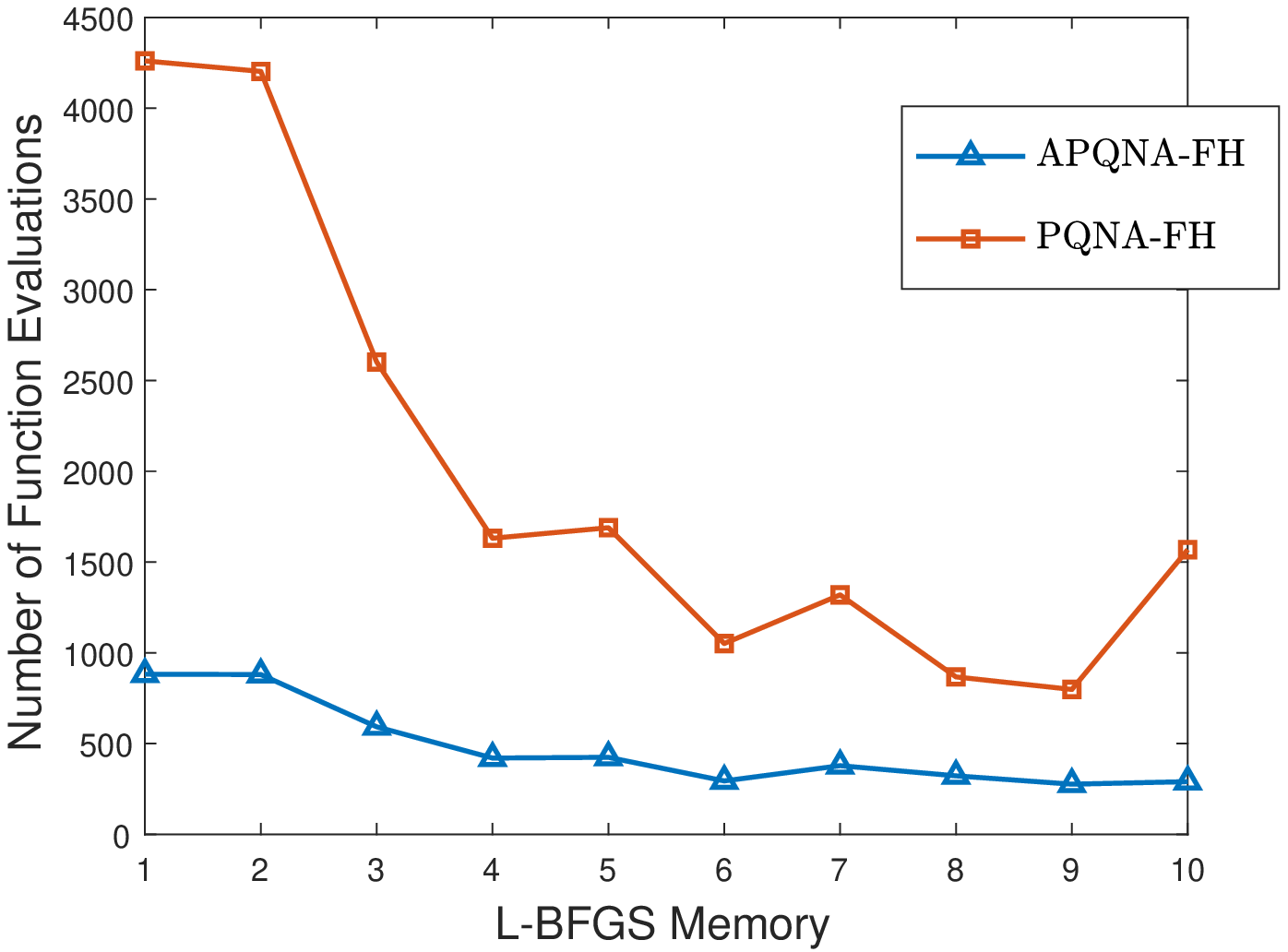}
  \caption{mnist (d=782, m=100000)}
\end{subfigure}
\begin{subfigure}{.51\textwidth}
  \centering
  \includegraphics[width=1\linewidth]{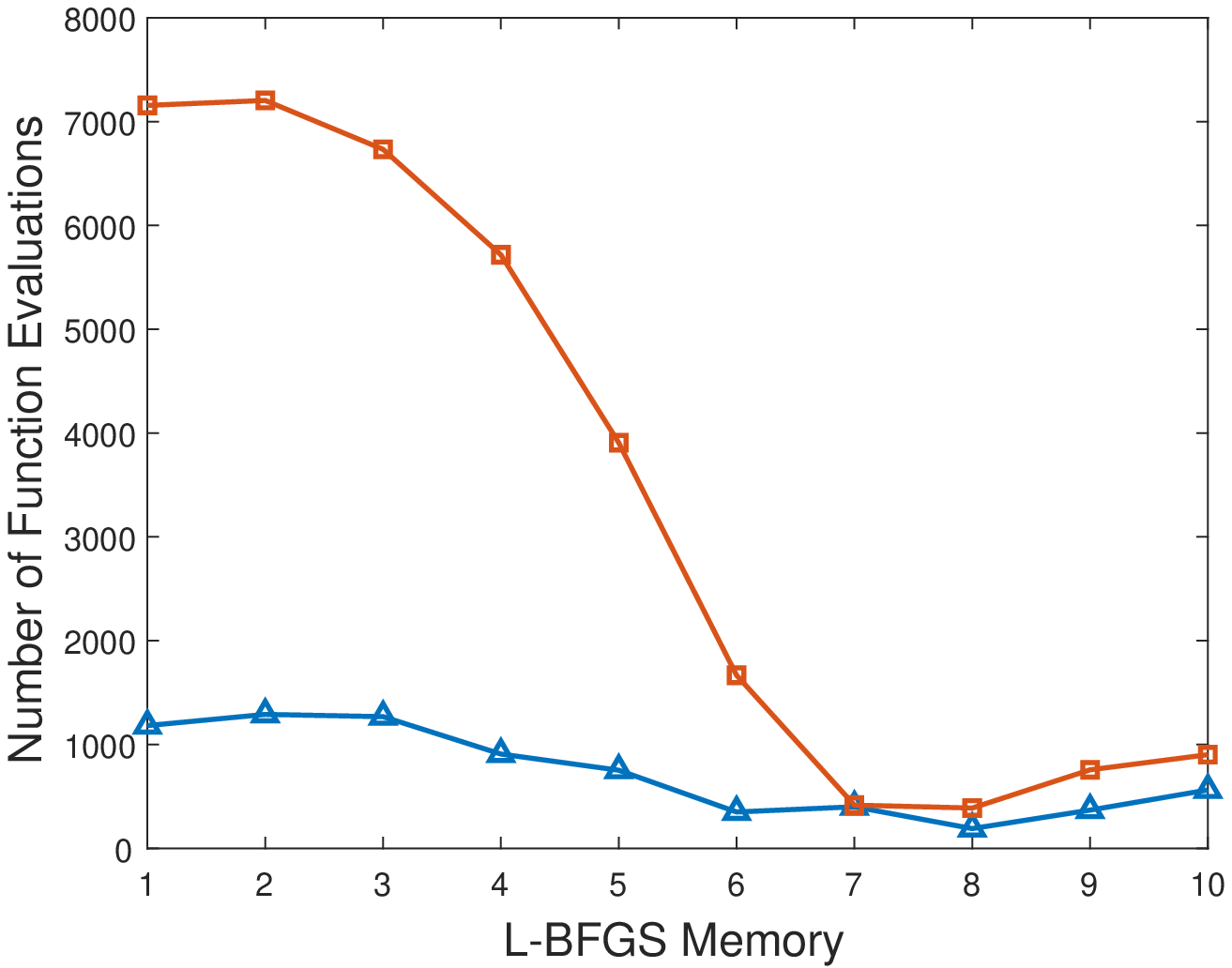}
  \caption{connect-4 (d=126, m=10000)}
\end{subfigure}
\begin{subfigure}{.51\textwidth}
  \centering
  \includegraphics[width=1\linewidth]{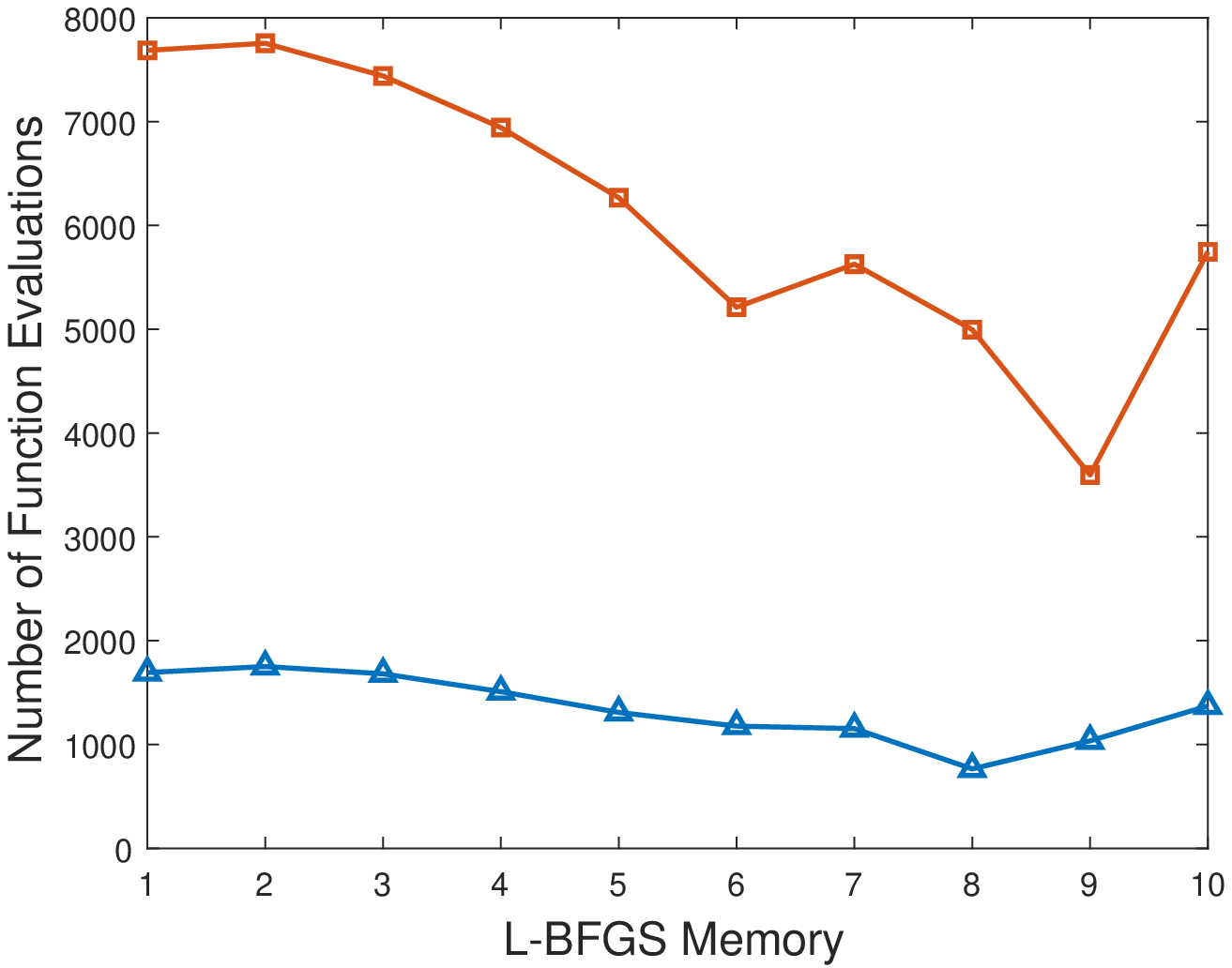}
  \caption{HAPT (d=561, m=7767)}
\end{subfigure}
\caption{{APQNA-FH} vs. {PQNA-FH} in terms of number of function evaluations.}
\label{p.fig3}
\end{figure}

The next experiment is to compare  the performance of {APQNA-FH} and {PQNA-FH} to observe the effect of acceleration in the fixed matrix setting. This comparison is done in terms of the number of iterations and the number of function evaluations, and is shown in Figure \ref{p.fig2} and Figure \ref{p.fig3}, respectively. The subproblem solution time is the same for both algorithms.
As we can see in Figure \ref{p.fig2}, in terms of the number of iterations, {APQNA-FH} dominates {PQNA-FH}, for a range of memory sizes of L-BFGS which have been used to compute matrix $H$. Moreover,  as is seen in Figure \ref{p.fig3}, {APQNA-FH} dominates {PQNA-FH}, in terms of the number of function evaluations, even though each iteration of {APQNA-FH} requires two function evaluations,  because of the nature of the accelerated scheme. This shows that {APQNA-FH} achieves practical acceleration compared to {PQNA-FH}, as supported by the theory in the previous section. 

Next, we compare the performance of {APQNA-FH} versus {APQNA-LBFGS} to compare the effect of  using the fixed approximate Hessian $H_k = \frac{1}{\sigma_k} H$, which satisfies condition $\sigma_{k} H_{k}  \preceq \sigma_{k-1} H_{k-1}$ versus using variable Hessian estimates computed via  L-BFGS method at each iteration,  while relaxing  condition $\sigma_{k} H_{k}  \preceq \sigma_{k-1} H_{k-1}$. Table \ref{p.table3} shows the results of this comparison, obtained based on the best choices of memory size for L-BFGS, in particular $\bar{k}=8$ and $\bar{k}=9$, respectively. As we can see, these two algorithms are competitive both in terms of the number of iterations and also the total solution time. Since {APQNA-FH} does not use the local information of function $f$ to approximate $H_k$, it often takes more iterations than {APQNA-LBFGS}, which constantly updates $H_k$ matrices. On the other hand, since {APQNA-FH} does not require additional computational effort to evaluate $H_k$, hence one iteration of this algorithm is  cheaper than one iteration of {APQNA-LBFGS}, which causes the competitive total solution time.

\btable[ht!]     
 \centering
 \captionsetup{justification=centering}
  \caption{{APQNA-FH} vs. {PQNA-LBFGS} in terms of function value (Fval), number of iterations (iter) and total solution time (time) in seconds.} 
  \label{p.table3}
 \small
  \begin{tabular}{ c|cc|cc|ccc}  
\multicolumn{8}{c}{\textbf{a9a}} \\ \hline
\abovespace
 Algorithm&iter&Fval& iter& Fval&iter&Fval&time\\ \hline
 \abovespace
 APQNA-LBFGS&20&3.4760e-01&40&3.4703e-01&64&3.4703e-01&2.83e+00\\ 
 \abovespace
APQNA-FH&20&3.4763e-01&40&3.4704e-01&99&3.4703e-01&4.33e+00\\ 
\multicolumn{8}{c}{\vspace{-0.1in}} \\
\multicolumn{8}{c}{\textbf{mnist}} \\ \hline
\abovespace
 Algorithm&iter&Fval& iter&Fval&iter&Fval&time\\ \hline
 \abovespace
 APQNA-LBFGS&50&8.9713e-02&100&8.9695e-02&148&8.9695e-02&1.04e+02\\ 
 \abovespace
APQNA-FH&50&8.9797e-02 &100&8.9698e-02&160&8.9695e-02&1.15e+02\\ 
\multicolumn{8}{c}{\vspace{-0.1in}} \\
\multicolumn{8}{c}{\textbf{connect-4}} \\ \hline
\abovespace
 Algorithm&iter&Fval& iter&Fval&iter&Fval&time\\ \hline
 \abovespace
APQNA-LBFGS&30&3.7769e-01&60&3.7688e-01&144&3.7682e-01&8.35e+00\\ 
\abovespace
APQNA-FH&30&3.7689e-01&60&3.7682e-01&93&3.7682e-01&3.95e+00\\ 
\multicolumn{8}{c}{\vspace{-0.1in}} \\
\multicolumn{8}{c}{\textbf{HAPT}} \\ \hline
\abovespace
 Algorithm&iter&Fval& iter& Fval&iter&Fval&time\\ \hline
 \abovespace
APQNA-LBFGS&120&7.1860e-02&240&7.1519e-02&356&7.1511e-02&1.04e+02\\ 
\abovespace
APQNA-FH&120&7.2134e-02&240&7.1523e-02&376&7.1511e-02&6.89e+01\\ 
\end{tabular} 
\centering
\etable

\begin{figure} [H]
\begin{subfigure}{.51\textwidth}
  \centering
  \includegraphics[width=1\linewidth]{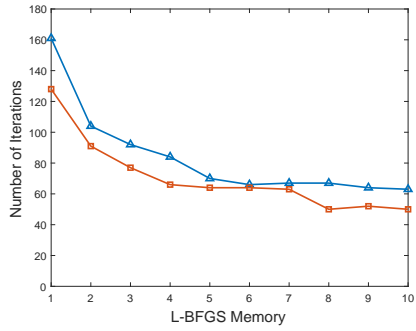}
  \caption{a9a (d=123, m=32561)}
  \label{fig:sfig1}
\end{subfigure}%
\begin{subfigure}{.51\textwidth}
  \centering
  \includegraphics[width=1\linewidth]{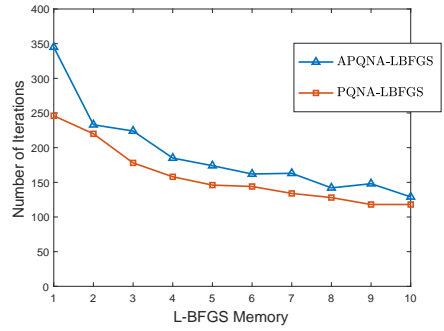}
  \caption{mnist (d=782, m=100000)}
  \label{fig:sfig2}
\end{subfigure}
\begin{subfigure}{.51\textwidth}
  \centering
  \includegraphics[width=1\linewidth]{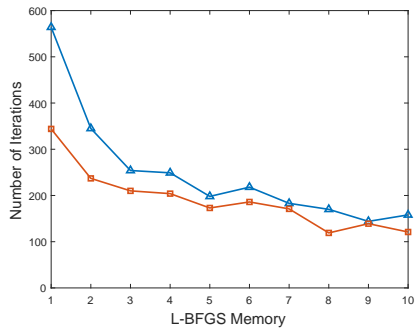}
  \caption{connect-4 (d=126, m=10000)}
  \label{fig:sfig2}
\end{subfigure}
\begin{subfigure}{.51\textwidth}
  \centering
  \includegraphics[width=1\linewidth]{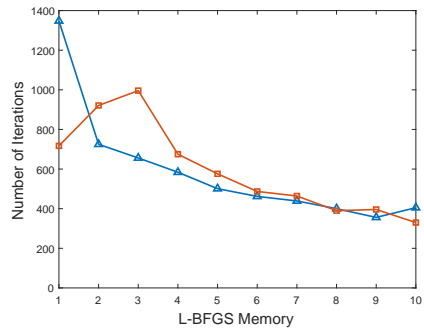}
  \caption{HAPT (d=561, m=7767)}
  \label{fig:sfig2}
\end{subfigure}
\caption{{APQNA-LBFGS} vs. {PQNA-LBFGS} in terms of number of iterations.}
\label{p.fig5}
\end{figure}

\begin{figure}[H]
\begin{subfigure}{.51\textwidth}
  \centering
  \includegraphics[width=1\linewidth]{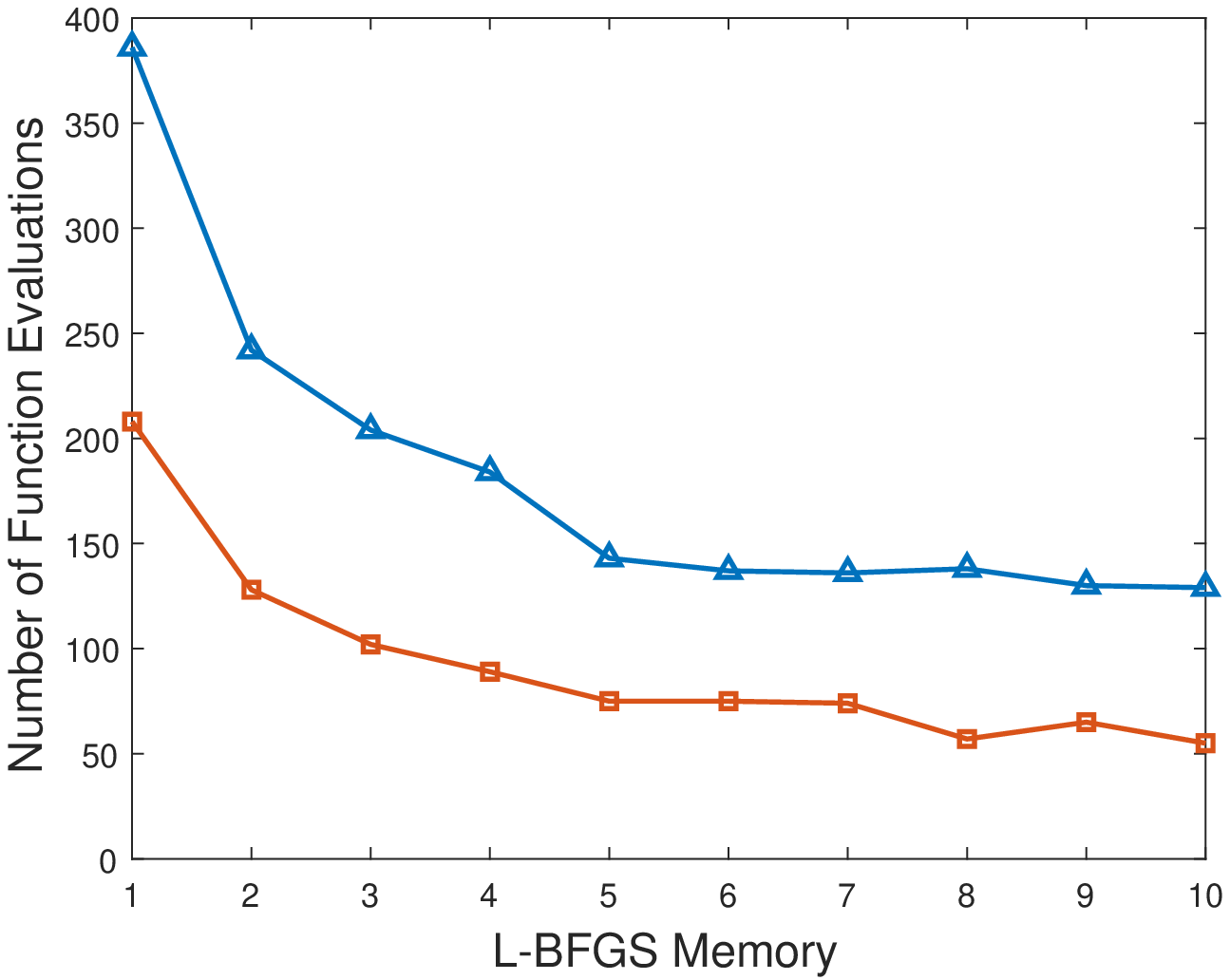}
  \caption{a9a (d=123, m=32561)}
  \label{fig:sfig1}
\end{subfigure}%
\begin{subfigure}{.51\textwidth}
  \centering
  \includegraphics[width=1\linewidth]{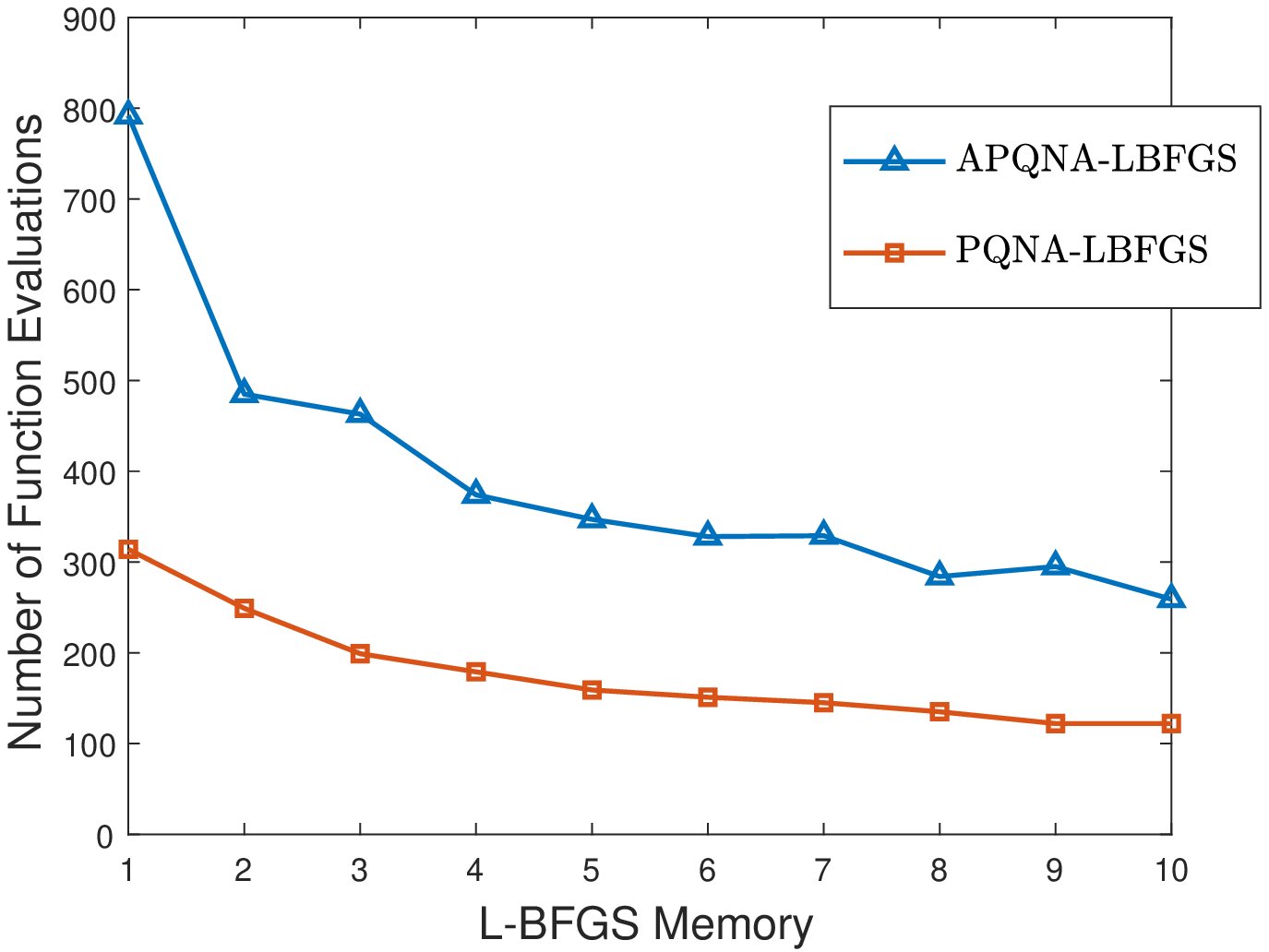}
  \caption{mnist (d=782, m=100000)}
  \label{fig:sfig2}
\end{subfigure}
\begin{subfigure}{.51\textwidth}
  \centering
  \includegraphics[width=1\linewidth]{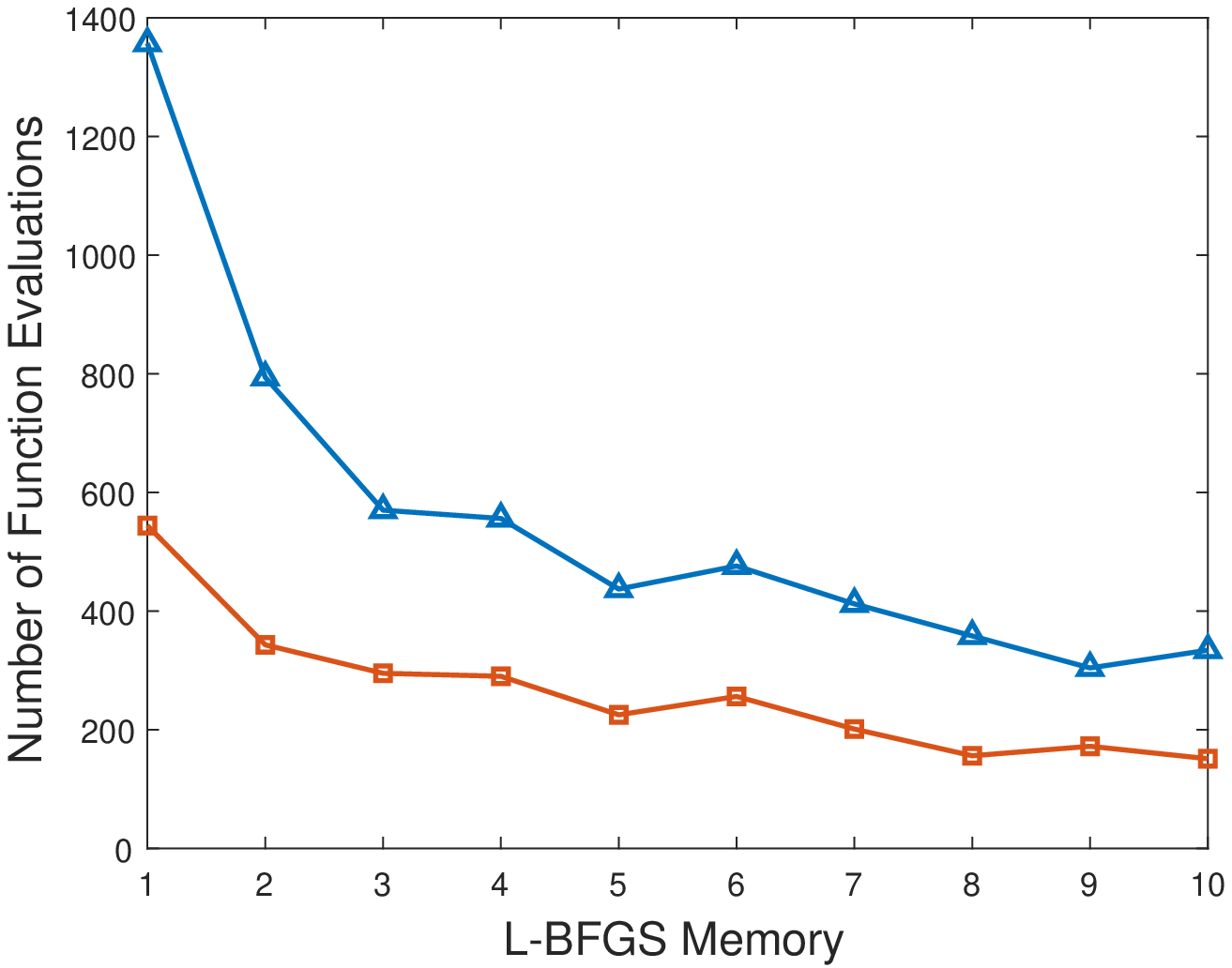}
  \caption{connect-4 (d=126, m=10000)}
  \label{fig:sfig2}
\end{subfigure}
\begin{subfigure}{.51\textwidth}
  \centering
  \includegraphics[width=1\linewidth]{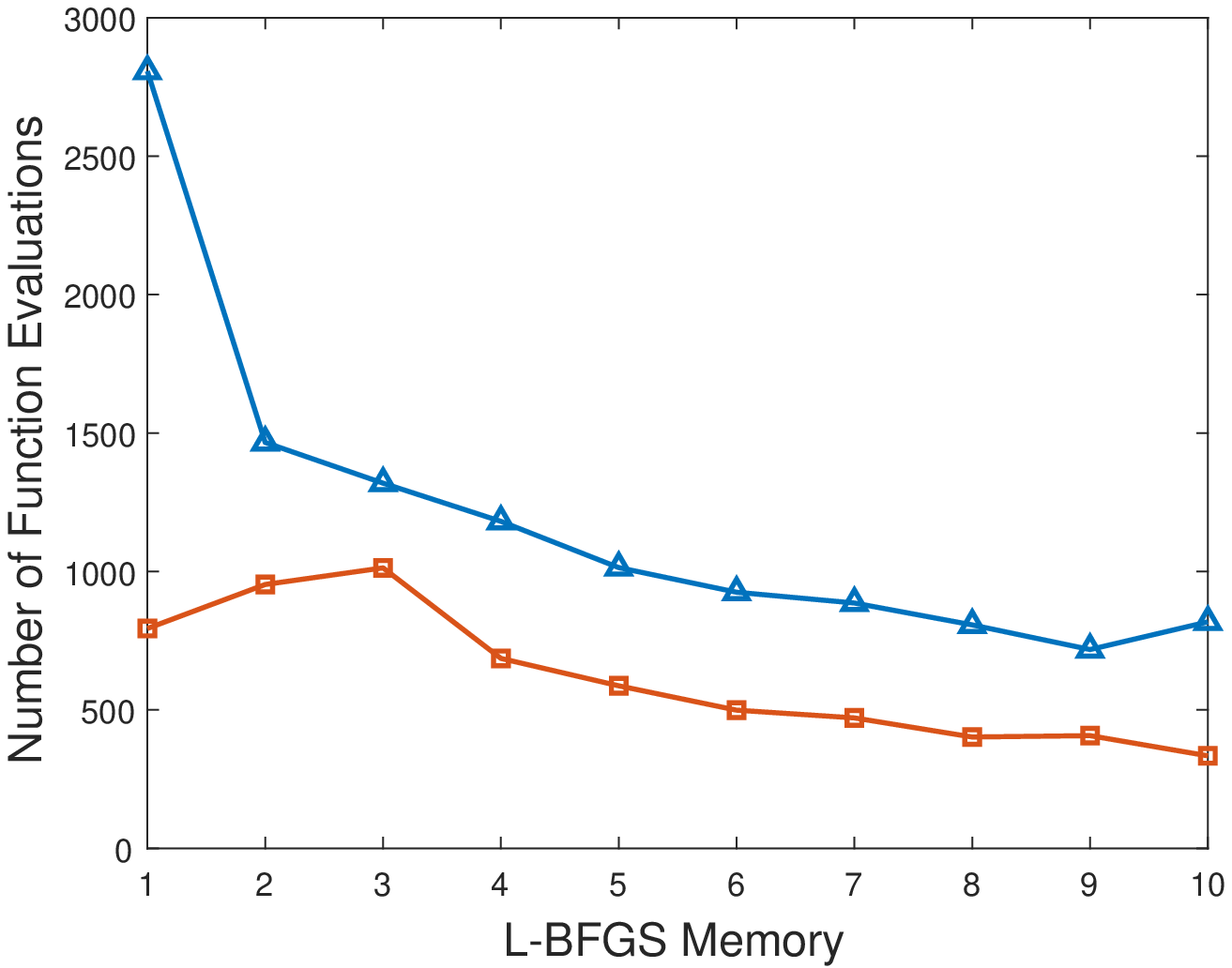}
  \caption{HAPT (d=561, m=7767)}
  \label{fig:sfig2}
\end{subfigure}
\caption{{APQNA-LBFGS} vs. {PQNA-LBFGS} in terms of number of function evaluations.}
\label{p.fig6}
\end{figure}

Finally, we compare {APQNA-LBFGS} and {PQNA-LBFGS}, to demonstrate the effect of using an accelerated scheme in the quasi-Newton type proximal algorithms. The results of this comparison are shown in Figure \ref{p.fig5} and Figure \ref{p.fig6} in terms of the number of iterations and the number of function evaluations,  respectively, for different memory sizes of L-BFGS Hessian approximation.

Clearly, as we can see in Figure \ref{p.fig5} and \ref{p.fig6}, not only does the accelerated scheme not achieve practical acceleration compared to {PQNA-LBFGS} in terms of the number of iterations, but it is also  inferior in terms of the number of function evaluations, since  every iteration requires two  function evaluations. Thus, we believe that the practical experiments support our theoretical analysis in that applying acceleration scheme in the case of  variable Hessian estimates may not result in a faster algorithm.  

\section{Conclusion}\label{sec.conclusion}
In this paper, we established a linear convergence rate of {PQNA} proposed in \cite{tang} under a strong convexity assumption.
To our knowledge, this is the first such result, for proximal quasi-Newton type methods, which have lately been popular in the literature. We also show that this convergence rate is preserved when subproblems are solved inexactly. We provide a simple and practical rule for the number of inner iterations that  guarantee sufficient accuracy of subproblem solutions. Moreover, we allow a relaxed  sufficient decrease condition during backtracking, which preserves the convergence rate, while it is known to improve the practical performance of the algorithm. 

Furthermore, we presented a variant of {APQNA} as an extension of {PQNA}.  We have shown that this algorithm has the convergence rate of $\mathcal{O}({1}/{k^2})$ under a strong condition on the Hessian estimates, which can not always be guaranteed in practice. We have shown that this condition holds
when  Hessian estimates are a multiple of a fixed matrix, which is computationally less expensive than the more common methods, such as the L-BFGS scheme. Although, this proposed algorithm has the same rate of convergence as the classic {APGA}, it is significantly faster in terms of the final number of iterations and also the total solution time. Based on the theory, using L-BFGS Hessian approximation, may result in worse convergence rate, however, our computational results show 
that the practical performance is about the same as that while the fixed matrix.  On the other hand, although in these two algorithms, we are applying the accelerated scheme, their practical performances are inferior to that of {PQNA-LBFGS}, which does not use any accelerated scheme and potentially has a slower sublinear rate of convergence in the absence of strong convexity. We conclude that using variable Hessian estimates is the most efficient approach, and will result in the linear convergence rate in the presence of strong convexity, but that a standard accelerated scheme is not useful in this setting. Exploring other possibly more effective accelerated schemes for the proximal quasi-Newton methods is the subject of future research.

\bibliographystyle{plain}
\bibliography{references}

\begin{thebibliography}{10}

\bibitem{beck}
A.~Beck and M.~Teboulle.
\newblock A fast iterative shrinkage-thresholding algorithm for linear inverse
  problems.
\newblock {\em SIAM}, 2:183--202, 2009.

\bibitem{nocedal3}
R.~Byrd, J.~Nocedal, and F.~Oztoprak.
\newblock An inexact successive quadratic approximation method for convex
  $\ell_1$-regularized optimization.
\newblock {\em Technical Report}, 2013.

\bibitem{nocedal1}
R.~H Byrd, J.~Nocedal, and R.~B Schnabel.
\newblock Representations of quasi \textsc{N}ewton matrices and their use in
  limited memory methods.
\newblock {\em Mathematical Programming}, 63:129--156, 1994.

\bibitem{dmitri}
D.~Drusvyatskiy and A.~S. Lewis.
\newblock Error bounds, quadratic growth, and linear convergence of proximal
  methods.
\newblock {\em Technical Report}, 2011.

\bibitem{phd}
H.~Ghanbari and K.~Scheinberg.
\newblock Optimization algorithms in machine learning.
\newblock {\em Doctoral Dissertation}, 2016.

\bibitem{hsieh}
C.~J. Hsieh, M.~Sustik, I.~Dhilon, and P.~Ravikumar.
\newblock Sparse inverse covariance matrix estimation using quadratic
  approximation.
\newblock {\em NIPS}, pages 2330--2338, 2011.

\bibitem{toh}
K.~Jiang, D.~Sun, and K.~Toh.
\newblock An inexact accelerated proximal gradient method for large scale
  linearly constrained convex \textsc{SDP}.
\newblock {\em SIAM}, 22:1042--1064, 2012.

\bibitem{lee}
J.~D. Lee, Y.~Sun, and M.~A. Saunders.
\newblock Proximal \textsc{N}ewton-type methods for convex optimization.
\newblock {\em Technical Report}, 2012.

\bibitem{yudin}
A.~Nemirovski and D.~Yudin.
\newblock Informational complexity and efficient methods for solution of convex
  extremal problems.
\newblock {\em J. Wiley and Sons, New York}, 1983.

\bibitem{nesterov3}
Y.~E. Nesterov.
\newblock A method for solving the convex programming problem with convergence
  rate $\mathcal{O} (1/k^2)$.
\newblock {\em Soviet Mathematics Doklady}, 27:372--376, 1983.

\bibitem{nesterov2}
Y.~E. Nesterov.
\newblock {\em Introductory Lectures on Convex Programming: A Basic Course}.
\newblock Springer, 2004.

\bibitem{nesterov5}
Y.~E. Nesterov.
\newblock Smooth minimization for non-smooth functions.
\newblock {\em Mathematical Programming}, 103:127--152, 2005.

\bibitem{nesterov4}
Y.~E. Nesterov.
\newblock Gradient methods for minimizing composite objective function.
\newblock {\em Mathematical Programming}, 140:125--161, 2013.

\bibitem{nocedal2}
J.~Nocedal and S.J. Wright.
\newblock {\em Numerical Optimization}.
\newblock Springer Series in Operations Research. Springer, New York, NY, USA,
  2nd edition, 2006.

\bibitem{olsen}
P.~A. Olsen, F.~Oztoprak, J.~Nocedal, and S.~J. Rennie.
\newblock Newton-like methods for sparse inverse covariance estimation.
\newblock {\em NIPS}, pages 764--772, 2012.

\bibitem{martin}
P.~Richtarik and M.~Takac.
\newblock Iteration complexity of randomized block-coordinate descent methods
  for minimizing a composite function.
\newblock {\em Mathematical Programming}, 144:1--38, 2014.

\bibitem{goldfarb}
K.~Scheinberg, D.~Goldfarb, and X.~Bai.
\newblock Fast first-order methods for composite convex optimization with
  backtracking.
\newblock {\em Foundation of Mathematics}, 14:389--417, 2014.

\bibitem{rish}
K.~Scheinberg and I.~Rish.
\newblock A greedy coordinate ascent method for sparse inverse covariance
  selection problem.
\newblock {\em SINCO, Technical Report}, 2009.

\bibitem{tang}
K.~Scheinberg and X.~Tang.
\newblock Practical inexact proximal quasi-\textsc{N}ewton method with global
  complexity analysis.
\newblock {\em Mathematical Programming}, 160:495--529, 2016.

\bibitem{mark}
M.~Schmidt, N.~L. Roux, and F.~Bach.
\newblock Convergence rate of inexact proximal-gradient method for convex
  optimization.
\newblock {\em NIPS}, pages 1458--1466, 2011.

\bibitem{bach}
M.~Schmidt, N.~L. Roux, and F.~Bach.
\newblock Supplementary material for the paper convergence rates of inexact
  proximal-gradient methods for convex optimization.
\newblock {\em NIPS}, 2011.

\bibitem{shalev}
S.~Shalev-Shwartz and A.~Tewari.
\newblock Stochastic methods for $\ell_1-$regularized loss minimization.
\newblock {\em ICML}, pages 929--936, 2009.

\bibitem{sra}
S.~Sra, S.~Nowozin, and S.J. Wright.
\newblock {\em Optimization for Machine Learning}.
\newblock Mit Pr, 2011.

\bibitem{tibshirani}
R.~Tibshirani.
\newblock Regression shrinkage and selection via the lasso.
\newblock {\em Journal of the Royal Statistical Society, Series B,
  Methodological}, 58:267--288, 1996.

\bibitem{tseng}
P.~Tseng.
\newblock On accelerated proximal gradient methods for convex-concave
  optimization.
\newblock {\em Technical Report}, 2008.

\bibitem{villa}
S.~Villa, S.~Salzo, L.~Baldassarre, and A.~Verri.
\newblock Accelerated and inexact forward-backward algorithms.
\newblock {\em SIAM}, 23:1607--1633, 2011.

\bibitem{lin}
G.~X. Yuan, K.~W. Chang, C.~J. Hsieh, and C.~J. Lin.
\newblock A comparison of optimization methods and software for large-scale
  $\ell_1$-regularized linear classification.
\newblock {\em JMLR}, 11:3183--3234, 2010.

\end{thebibliography}
\end{document}